\PassOptionsToPackage{lowtilde}{url} 
\documentclass[a4paper]{siamart}

\usepackage[utf8]{inputenx}      
\usepackage[T1]{fontenc}         
\usepackage{grffile}             
\usepackage{euscript,graphicx,epsfig,amsfonts,mathtools,epstopdf}
\usepackage[ruled,lined,vlined,linesnumbered,algo2e]{algorithm2e}
\usepackage{siunitx}
\usepackage{listings}

\usepackage{xcolor}
\definecolor{darkgreen}{rgb}{0.0,0.6,0.0}
\usepackage[draft]{changes}

\def\cal{\mathcal}
\def\abs#1{\lvert#1\rvert}
\def\norm#1{\left\|#1\right\|}
\providecommand{\bignorm}[1]{\big\lVert#1\big\rVert}

\def\R{\mathbb{R}}

\def\Rn{\R^n}

\def\Rmm{\R^{m\times m}}
\def\Rnn{\R^{n\times n}}

\def\BA{\textbf{A}}  \def\CA{{\cal A}}

  \def\CK{{\cal K}}

  \def\CP{{\cal P}}

\def\Bzero{\boldsymbol{0}}
\def\BAs{\BA{\kern-1.5pt}}

\def\CPs{\CP{\kern-0.8pt}}

\mathcode`\@="8000
{\catcode`\@=\active \gdef@{\mkern1mu}}

\def\mydate{\number\day\ {\ifcase\month \or January\or February\or
              March\or April\or May\or June\or July\or August\or
              September\or October\or November\or December\fi}
\number\year}
\def\vek#1{\mathbf{#1}}

\def\ip#1{\left\langle#1\right\rangle}
\providecommand{\bigip}[1]{\big\langle#1\big\rangle}
\providecommand{\argmin}[1]{\underset{#1}{\text{{\rm argmin}}}}

\providecommand{\prn}[1]{\left(#1\right)}
\providecommand{\bigprn}[1]{\big(#1\big)}
\providecommand{\Bigprn}[1]{\Big(#1\Big)}

\def\curl#1{\left\{#1\right\}}
\providecommand{\ab}[1]{\left|#1\right|}

\newcommand\restr[2]{{
  \left.\kern-\nulldelimiterspace 
  #1 
  \vphantom{\big|} 
  \right|_{#2} 
  }}

\def\Span{\rm span}

\def\be{\begin{equation}}
\def\ee{\end{equation}}
\def\bea{\begin{eqnarray}}
\def\eea{\end{eqnarray}}
\def\nn{\nonumber}
\def\mand{\mbox{\ \ \ and\ \ \ }}

\def\mwith{\mbox{\ \ \ with\ \ \ }}
\def\mwhere{\mbox{\ \ \ where\ \ \ }}

\def\msuchthat{\mbox{\ \ \ such that\ \ }}

\def\bbmat{\begin{bmatrix}}
\def\ebmat{\end{bmatrix}}

\def\balg{\begin{algorithm2e}\DontPrintSemicolon}
\def\ealg{\end{algorithm2e}}
\def\bthm{\begin{theorem}}
\def\ethm{\end{theorem}}
\def\blem{\begin{lemma}}
\def\elem{\end{lemma}}
\def\bprop{\begin{proposition}}
\def\eprop{\end{proposition}}
\def\bcor{\begin{corollary}}
\def\ecor{\end{corollary}}
\def\bdefin{\begin{definition}}
\def\edefin{\end{definition}}
\def\bc{\begin{cases}}
\def\ec{\end{cases}}
\def\bproof{\par\addvspace{1ex} \indent\textit{Proof}.\ \ }
\def\eproof{\hfill\cvd\linebreak\indent}
\newtheorem{exple}{Example}
\def\bex{\begin{exple}}
\def\eex{\end{exple}}
\def\bass{\begin{assumption}}
\def\eass{\end{assumption}}

\def\cvd{~\vbox{\hrule\hbox{%
  \vrule height1.3ex\hskip0.8ex\vrule}\hrule } }

\providecommand{\dx}{\textup{d}\boldsymbol{x}}

\DeclareMathSymbol{\dprod}{\mathbin}{operators}{"3A}
\providecommand{\blkdiag}{\operatorname{blkdiag}}
\providecommand{\python}{Python}
\providecommand{\fenics}{\textsc{FEniCS}}
\providecommand{\petsc}{\textsc{PETSc}}
\providecommand{\matlab}{\textsc{Matlab}}
\providecommand{\file}[1]{\texttt{\nolinkurl{#1}}}
\providecommand{\laplace}{\triangle}

\newsiamthm{example}{Example}

\renewcommand{\theequation}{\thesection.\arabic{equation}} 
\renewcommand{\thefigure}{\thesection.\arabic{figure}}
\renewcommand{\thetable}{\thesection.\arabic{table}}
\crefname{algocf}{algorithm}{algorithms}
\Crefname{algocf}{Algorithm}{Algorithms}
\numberwithin{theorem}{section}

\newcommand{\TheTitle}{A modified implementation of MINRES to monitor residual subvector norms for block systems} 
\newcommand{\TheAuthors}{R. Herzog and K. M. Soodhalter}

\headers{MONITORING OF RESIDUAL IN SADDLE-POINT PROBLEMS}{\TheAuthors}

\title{{\TheTitle}\thanks{This version dated \today.}}

\author{
	Roland Herzog\thanks{Technische Universität Chemnitz, Faculty of Mathematics, Professorship Numerical Mathematics (Partial Differential Equations), D--09107 Chemnitz, Germany (\email{roland.herzog@mathematik.tu-chemnitz.de}, \url{https://www.tu-chemnitz.de/herzog}).}
	\and
	Kirk M. Soodhalter\thanks{Transfer Group, Johann Radon Institute for Computational and Applied Mathematics, Altenbergerstraße~69, A--4040 Linz, Austria (\email{kirk.soodhalter@ricam.oeaw.ac.at}, \url{http://math.soodhalter.com}).}
}

\ifpdf
\hypersetup{
	pdftitle={\TheTitle},
	pdfauthor={\TheAuthors}
}
\fi

\begin{document}

\maketitle


\begin{abstract}
	Saddle-point systems, i.e., structured linear systems with symmetric matrices are considered. 
	A modified implementation of (preconditioned) MINRES is derived which allows to monitor the norms of the subvectors individually.
		Compared to the implementation from the textbook of [Elman, Sylvester and Wathen, Oxford University Press, 2014], our method requires one extra vector of storage and no additional applications of the preconditioner.
		Numerical experiments are included.
\end{abstract}



\section{Introduction}\label{section.intro}
We are solving symmetric linear systems of the form
\be\label{eqn.ABupf}
	\bbmat \vek A & \vek B^{T}\\ \vek B & -\vek C \ebmat \bbmat \vek u \\ \vek p \ebmat = \bbmat \vek f_{\vek u}\\ \vek f_{\vek p} \ebmat,
\ee
where $\vek A\in\Rmm$ and $\vek C\in\R^{p\times p}$ are symmetric, and $\vek B\in\R^{p\times m}$,
by applying MINRES or preconditioned MINRES \cite{Paige1975} iteration.
 After $j$ iterations, we have approximation $(\vek u^{(j)},\vek p^{(j)})$ and residual 
\be\nn
	\vek r^{(j)} = \bbmat \vek f_{\vek u}\\ \vek f_{\vek p} \ebmat - \bbmat \vek A & \vek B^{T}\\ \vek B & -\vek C \ebmat \bbmat \vek u^{(j)} \\ \vek p^{(j)} \ebmat = \bbmat \vek f_{\vek u}\\ \vek f_{\vek p} \ebmat - \bbmat \vek A\vek u^{(j)} + \vek B^{T}\vek p^{(j)}\\ \vek B\vek u^{(j)}  - \vek C\vek p^{(j)}\ebmat.
\ee
We denote the two parts of the residual as
\be\nn
	\vek r^{(j)}_{\vek u} =  \vek f_{\vek u} - \bigprn{\vek A\vek u^{(j)} + \vek B^{T}\vek p^{(j)}}\mand \vek r^{(j)}_{\vek p} = \vek f_{\vek p} - \bigprn{\vek B\vek u^{(j)}  - \vek C\vek p^{(j)}}.
\ee
The question we seek to answer is: can we monitor the individual preconditioned norms of $\vek r^{(j)}_{\vek u}$ and $\vek r^{(j)}_{\vek p}$ (as opposed to the full norm of $\vek r^{(j)}$)
using only quantities arising in an efficient implementation of
the preconditioned MINRES algorithm, namely that in \cite[Algorithm~4.1]{ESW-Book.2014}?  The answer is yes, under certain
conditions.  The technical requirements to do so are connected to the notion of so-called ``natural'' preconditioners
	which arise in this setting; see, e.g., \cite{Pestana.Wathen.Natural-precond.2015,Zulehner.Nonstandard-norms.2011} and references contained therein. In \cref{section.deriv}, we demonstrate that at the storage cost of one additional
full-length vector and six scalars but no additional applications of the preconditioner, one can modify the preconditioned MINRES method to calculate these norms in a progressive fashion.
	An extension to more than two residual parts is straightforward.
	An implementation of our modified version of MINRES is given in \cref{alg.prec.minres-modified} and is available at \cite{SUBMINRES} as a Matlab file.

	As a motivation for our study, we mention that the individual parts of the residual in \eqref{eqn.ABupf} often have different physical interpretations.
	Monitoring them individually allows a better insight into the convergence of MINRES and it allows the formulation of refined stopping criteria. 
	We present examples and numerical experiments in \cref{section.examples}.

\section{How to monitor both parts of the preconditioned residual}\label{section.deriv}
In this section, we derive how one monitors these norms without incurring much
extra computational or storage expense.  
We will describe everything in terms of the variable names use in \cite[Algorithm~4.1]{ESW-Book.2014} with two exceptions,
which will be noted below.

Here we describe quickly the derivation of preconditioned MINRES as a Krylov subspace method, specifically in order to relate certain quantities
from \cite[Algorithm~4.1]{ESW-Book.2014} to the common quantities arising in the Lanczos process.
In principle, the preconditioned MINRES algorithm is an implementation of the minimum residual Krylov subspace method for 
matrices which are symmetric.  It can be equivalently formulated as an iteration for operator
equations in which the operator is self-adjoint, mapping elements of a Hilbert space into its dual; see \cite{G.Herzog.S.2014}.  We present the derivation for the finite dimensional case.
For the purposes of this discussion, let $\vek x^{(j)} = (\vek u^{(j)}, \vek p^{(j)})$ be the $j$th approximation.  
Let $\vek K = \bbmat \vek A & \vek B^{T}\\ \vek B & \vek 0 \ebmat$. 
We further assume that the preconditioner $\vek P$ is symmetric positive definite and has block diagonal structure
\be\nn
	\vek P = \bbmat \vek P_{\vek u} & \\ & \vek P_{\vek p} \ebmat \mwith \vek P_{\vek u}\in\R^{m\times m} \mand \vek P_{\vek p}\in\R^{p\times p}.
\ee
This assumption is natural in many situations when we acknowledge that the residual subvectors often have different physical meaning and need to be measured in different norms; see \cite{Pestana.Wathen.Natural-precond.2015,Zulehner.Nonstandard-norms.2011} and \cref{section.examples} for examples.

We briefly review some basic facts about Krylov subspaces for symmetric matrices.  
We begin with the unpreconditioned situation.
Let $\boldsymbol \CA \in \Rnn$ be symmetric.  For some starting element $\vek h\in\Rn$, we define the $j$th Krylov subspace generated by
$\boldsymbol\CA$ and $\vek h$ to be 
\be\nn
	\CK_{j}(\boldsymbol\CA, \vek h) = \Span\curl{\vek h, \boldsymbol\CA \vek h, \boldsymbol\CA^{2}\vek h, \ldots, \boldsymbol\CA^{j-1}\vek h}
\ee
using the Lanczos process, where if $\vek V_{j} = \bbmat \vek v_{1} & \vek v_{2} & \cdots & \vek v_{j} \ebmat\in\R^{n\times j}$
is the matrix with columns that form an orthonormal basis of $\CK_{j}(\boldsymbol\CA, \vek h) $, then we have the Lanczos relation
\be\label{eqn.lanczos-rel}
	\boldsymbol\CA\vek V_{j} = \vek V_{j+1}\overline{\vek T}_{j}\mwhere \overline{\vek T}_{j}\in\R^{(j+1)\times j}.
\ee
The matrix $\overline{\vek T}_{j}$ is tridiagonal, and the matrix $\vek T_{j}$ (defined as the first $j$ rows of $\overline{\vek T}_{j}$) is
symmetric. 
This implies that to generate each new Lanczos vector, we must orthogonalize the newest vector against only the two 
previous Lanczos vectors.  This leads to the following well-known three-term recurrence formula
\be\label{eqn-3-term}
	\boldsymbol\CA\vek v_{j} = \gamma_{j+1}\vek v_{j+1} + \delta_{j}\vek v_{j} + \gamma_{j}\vek v_{j-1}.
\ee
Using the naming conventions in \cite[Algorithm~2.4]{ESW-Book.2014}, we have 
\be\nn
	\overline{\vek T}_{j} = \bbmat \delta_{1} & \gamma_{2} & &\\ \gamma_{2}& \ddots & \ddots & \\ & \ddots & \ddots & \gamma_{j} \\ & & \gamma_{j} & \delta_{j} \\ & & & \gamma_{j+1} \ebmat.
\ee

In the case that we have a symmetric positive definite preconditioner $\vek P$, we show that a Krylov subspace can still constructed using
the short-term Lanczos iteration and that a MINRES method can be used for solving \eqref{eqn.ABupf}.  In, e.g., \cite{ESW-Book.2014}, preconditioned 
MINRES is derived by first observing that $\vek P$ admits the Cholesky decomposition $\vek P = \vek H\vek H^{T}$.  Thus one could consider 
solving the two-sided preconditioned equations
\be\nn
	\vek H^{-1}\vek K \vek H^{-T}\vek y = \vek H^{-1}\vek b\mwith\vek y = \vek H^{T}\vek x \mwhere\vek x =\bbmat \vek u \\ \vek p \ebmat \mand\vek b =\bbmat \vek f_{1}\\ \vek f_{2} \ebmat,
\ee
where we have the initial approximation $\vek y^{(0)} = \vek H^{T}\vek x^{(0)}$.
 The matrix $\vek H^{-1}\vek K \vek H^{-T}$ is still symmetric.
However, one can notice that the preconditioned residual satisfies  
$\vek H^{-1}\vek r^{(k)} = \norm{\vek r^{(k)}}_{\vek P^{-1}}$ where $\norm{\cdot}_{\vek P^{-1}}$ is the norm arising from the inner product
induced by $\vek P^{-1}$, i.e., $\ip{\cdot,\cdot}_{\vek P^{-1}}=\ip{\vek P^{-1}\cdot,\cdot}$.  One further notes that we have the equivalence
\be\nn
	\CK_{j}(\vek H^{-1}\vek K \vek H^{-T},\vek H^{-1}\vek r^{(0)}) = \vek H^{-1}\CK_{j}(\vek K\vek P^{-1},\vek r^{(0)}).
\ee
If preconditioned MINRES at iteration $j$ produce a correction 
\be\nn
	\vek s^{(j)}\in\CK_{j}(\vek H^{-1}\vek K \vek H^{-T},\vek H^{-1}\vek r^{(0)})\msuchthat \vek y^{(j)} = \vek y^{(0)} + \vek s^{(j)}, 
\ee	
then we recover the
approximation for $\vek x$ with $\vek x^{(j)} = \vek H^{-T}\prn{\vek y^{(0)}+ \vek s^{(j)}}$.  
Since $\vek x^{(0)} = \vek H^{-T}\vek y^{(0)}$ the correction space for preconditioned MINRES with respect to the original variables 
is actually $\vek P^{-1}\CK_{j}(\vek K\vek P^{-1},\vek r^{(0)})$.  From here, one can
show that the preconditioned MINRES iteration is equivalent to an iteration with the subspace 
$\CK_{j}(\vek K\vek P^{-1},\vek r^{(0)})$ but with respect to the
inner product $\ip{\cdot,\cdot}_{\vek P^{-1}}$ induced by the preconditioner.  
This is indeed a short-term recurrence iteration, as it can
be shown that $\vek K\vek P^{-1}$ is symmetric with respect to $\ip{\cdot,\cdot}_{\vek P^{-1}}$.  Thus the derivation of the preconditioned
MINRES method \cite[Algorithm~4.1]{ESW-Book.2014} can be presented as a small modification of the unpreconditioned 
MINRES \cite[Algorithm~2.4]{ESW-Book.2014}.

This same fact was also shown in \cite{G.Herzog.S.2014} but by interpreting the underlying operator as a mapping from a Hilbert space to
its dual.  
There are many benefits to this interpretation. 
One advantage is that even in the finite dimensional situation, it allows the derivation of preconditioned MINRES without the temporary introduction of Cholesky factors.

Let the columns of $\vek V_{j}$ still form an orthonormal basis for a Krylov subspace, but this time the space 
$\CK_{j}(\vek K\vek P^{-1},\vek r^{(0)})$, and the orthonormality is with respect to $\ip{\cdot,\cdot}_{\vek P^{-1}}$.  
Let $\vek Z_{j} = \vek P^{-1}\vek V_{j}$ have as columns the image of those vectors
under the action of the preconditioner.  Then we have the preconditioned Lanczos relation,
\be\label{eqn.Lanczos-rel-prec}
	\vek K\vek P^{-1}\vek V_{j} = \vek K\vek Z_{j} = \vek V_{j}\overline{\vek T}_{j}
\ee
where $\overline{\vek T}_{j}$ has the tridiagonal structure described earlier.  What differs here is that the columns of 
$\vek V_{j}$ have been orthonormalized with respect to the $\vek P^{-1}$ inner product, and the entries of the tridiagonal
matrix $\overline{\vek T}_{j}$ are formed from these $\vek P^{-1}$ inner products.  We have the three-term recurrence
\be\label{eqn.3-term-rec-prec}
	\vek K\vek z_{j} = \gamma_{j+1}\vek v_{j+1}+\delta_{j}\vek v_{j} + \gamma_{j}\vek v_{j-1}.
\ee
The preconditioned MINRES method solves the following residual minimization problem, 
\be\label{eqn.min-res}
\vek x^{(j)} = \vek x^{(0)} + \vek Z_{j}\vek y^{(j)}\mwhere \vek y^{(j)} = \argmin{\vek y\in\R^{j}}\norm{\gamma_{1}\vek e_{1}-\overline{\vek T}_{j}\vek y}_{2},
\ee
where $\gamma_{1} = \norm{\vek r^{0}}_{\vek P^{-1}}$.  

Though this reduces the minimization of the residual to a small least-squares problem, this is not the most efficient way to implement 
the MINRES method; and indeed, this is not what is done in, e.g., \cite[Algorithm~4.1]{ESW-Book.2014}.  
Due to \eqref{eqn-3-term}, MINRES can be derived such that only six full-length vectors must be stored.
In the interest of not again deriving everything, we will 
simply provide a few relationships between quantities from the above explanation of MINRES and those in \cite[Algorithm~4.1]{ESW-Book.2014}.
Further implementation details can be found in, e.g., \cite[Section~2.5]{Greenbaum1997a}.

Let $\overline{\vek T}_{j} = \vek Q_{j}\overline{\vek R}_{j}$ be a QR-factorization of 
$\overline{\vek T}_{j}$ computed with Givens rotations where $\vek Q_{j}\in\R^{(j+1)\times (j+1)}$
is a unitary matrix constructed from the product of Givens rotations, and $\overline{\vek R}_{j}\in\R^{(j+1)\times j}$ is upper triangular.  The matrix $\vek R_{j}\in\R^{j\times j}$ is simply
$\overline{\vek R}_{j}$ with the last row (of zeros) deleted.  We can write the matrix 
$\vek Q_{j}^{T}$ as the product of the Givens rotation used to triangularize $\overline{\vek T}_{j}$.
Since $\overline{\vek T}_{j}$ is tridiagonal and Hessenberg, we must annihilate only one subdiagonal entry
per column.  Following from \cite[Algorithm~4.1]{ESW-Book.2014}, we denote $s_{i}$ and $c_{i}$
to be the Givens sine and cosine used to annihilate the $i$th entry of column $i-1$ using the
unitary matrix $\vek F_{i} = \bbmat c_{i} & s_{i} \\ -s_{i} & c_{i} \ebmat$.  Thus we can denote
\be\nn
	\vek Q_{j}^{T} = \vek G_{j+1}^{(j)}\vek G_{j}^{(j)}\cdots \vek G_{2}^{(j)}
\ee
where we define $\vek G_{i}^{(j)}\in\R^{(j+1)\times (j+1)}$ to be the matrix applying the $i$th 
Givens rotation to rows $i-1$ and $i$ to a $(j+1)\times (j+1)$ matrix.  
Using a normal equations formulation of \eqref{eqn.min-res},
we can derive an expression for the least squares minimizer 
\be\label{eqn.y-expression}
	\vek y^{(j)} = \vek R_{j}^{-1}\curl{\vek Q_{j}^{T}(\gamma_{1}\vek e_{1})}_{1:j},
\ee
where $\curl{\cdot}_{1:j}$ indicates we take only the first $j$ rows of the argument.
For the purpose of discussion, we make an additional modification to the variable naming using in 
\cite[Algorithm~4.1]{ESW-Book.2014}.  We index $\eta$, which is used to track the residual norm.  
On Line~4 of the algorithm, we would have $\eta_{0} = \gamma_{1}$, and at Line~18, 
we would write $\eta_{j} = -s_{j+1}\eta_{j-1}$.  From the proof of 
\cite[Corollary 2.5.3]{Fischer1996}, one can see that
at iteration $j$, the least squares residual can be written
\be\label{eqn.ls-resid}
	\gamma_{1}\vek e_{1} - \overline{\vek T}_{j}\vek y^{(j)} = \vek Q_{j}(\eta_{j}\vek e_{j+1}),
\ee
where $\vek e_{j+1}$ is the last column of the $(j+1)\times (j+1)$ identity matrix, and because
$\vek Q_{j}$ is unitary, it follows that $\ab{\eta_{j}}$ is the norm of the residual $\bignorm{\vek r^{(j)}}_{\vek P^{-1}}$ pertaining to $\vek x^{(j)}$.

We now derive an expression for the norms of $\vek r^{(j)}_{\vek u}$ and $\vek r^{(j)}_{\vek p}$ in terms of quantities arising in the
Lanczos process and residual minimization.  Due to the block structure of $\vek K$, we partition 
$\vek V_{j}$ and $\vek Z_{j}$ similarly,
\be\label{eqn.laczos-vec-part}
	\vek V_{j} = \bbmat \vek V_{j,\vek u}\\ \vek V_{j,\vek p}  \ebmat \mand \vek Z_{j} = \bbmat \vek Z_{j,\vek u}\\ \vek Z_{j,\vek p}\ebmat
\ee
and then insert the partitioned vectors into \eqref{eqn.lanczos-rel}
\be\nn
	\bbmat \vek A & \vek B^{T}\\ \vek B & -\vek C \ebmat \bbmat \vek Z_{j,\vek u}\\ \vek Z_{j,\vek p}  \ebmat = \bbmat  \vek A\vek Z_{j,\vek u} + \vek B^{T}\vek Z_{j,\vek p}\\ \vek B\vek Z_{j,\vek u} - \vek C\vek Z_{j,\vek p} \ebmat = \bbmat \vek V_{j+1,\vek u} \\ \vek V_{j+1,\vek p} \ebmat \overline{\vek T}_{j}= \bbmat \vek V_{j+1,\vek u}\overline{\vek T}_{j} \\ \vek V_{j+1,\vek p}\overline{\vek T}_{j} \ebmat .
\ee
We note that for tracking and updating the full residual vector, this partition is artificially imposed. One could simply track
the full residual and partition at the end to compute the subvector norms separately.  We derive everything using this block partition because
one may wish to track the norm of only one subvector.  Furthermore, one needs the partition quantities to define the recursions for 
updating the preconditioned subvector norms.
It then follows that we have Lanczos relations for the blocks
\be\label{eqn.part-lanczos-rel}
\vek A\vek Z_{j,\vek u} + \vek B^{T}\vek Z_{j,\vek p} = \vek V_{j+1,\vek u}\overline{\vek T}_{j}  \mand \vek B\vek Z_{j,\vek u} - \vek C\vek Z_{j,\vek p}  =\vek V_{j+1,\vek p}\overline{\vek T}_{j}.
\ee
From \cite[Theorem 2.5.7]{Fischer1996}, we know that the $j$th MINRES residual satisfies the recursion
	\be\label{eqn.minres-recursion}
	\vek r^{(j)} = \vek r^{(j-1)} - c_{j+1}\eta_{j-1}\vek m_{j}^{(j+1)},
	\ee
	where the auxiliary vectors can also be progressively updated.  The progressive updating scheme for the auxiliary vectors can also be seen
	by observing that these auxiliary vectors form an orthonormal basis for the constraint space 
	$\vek K\vek P^{-1}\CK_{j}(\vek K\vek P^{-1},\vek r^{(0)})$, namely
	\be\nn
		\vek M_{j+1} = \vek V_{j+1}\vek Q_{j} = \bbmat \vek m_{1}^{(j+1)} & \vek m_{2}^{(j+1)} & \cdots & \vek m_{j}^{(j+1)} & {\vek m}_{j+1}^{(j+1)} \ebmat.
	\ee
Using induction, one can show that since $\vek Q_{j}$ is the product of Givens rotations  with $\vek Q_{j} = \vek Q_{j-1}\vek G_{j+1}^{(j)T}$, we have that 
$\vek m_{i}^{(j+1)} = \vek m_{i}^{(j)}$ for all $1 \leq i\leq j-1$, and 
\be\nn
	\vek m_{j}^{(j+1)} = c_{j+1}{\vek m}_{j}^{(j)} + s_{j+1}\vek v_{j+1}\mand {\vek m}_{j+1}^{(j+1)} = -s_{j+1}{\vek m}_{j}^{(j)} + c_{j+1}\vek v_{j+1},
\ee
where we set $\vek m_{1}^{(1)} = \vek v_{1}$.	
This implies that we can drop the superindices for the first $j-1$ columns of $\vek M_{j+1}$ 
which remain unchanged, i.e.,
\bea
	\vek M_{j+1} &=& \bbmat \vek m_{1} & \vek m_{2} & \cdots & \vek m_{j-1} & \vek m_{j}^{(j+1)} & {\vek m}_{j+1}^{(j+1)} \ebmat.
\eea
From \eqref{eqn.ls-resid}, we then have that
\be\nn
	\vek r^{(j)} = \eta_{j}\vek V_{j+1}\vek Q_{j}\vek e_{j+1} = \eta_{j}{\vek m}_{j+1}^{(j+1)}.
\ee
If we then partition the vector $\vek{m}_{j+1}^{(j+1)} = \bbmat \vek m_{j+1,\vek u}^{(j+1)} \\ \overset{}{\vek m_{j+1,\vek p}^{(j+1)}}\ebmat$, we then have formulas for the partial residuals,
\be\label{eqn.Buj-direct}
		\vek r_{\vek u}^{(j)} = \eta_{j}\vek m_{j+1,\vek u}^{(j+1)} \mand \vek r_{\vek p}^{(j)} = \eta_{j}\vek m_{j+1,\vek p}^{(j+1)}.
\ee
%
%
We can also rewrite \eqref{eqn.Buj-direct} in order to recover \eqref{eqn.minres-recursion}. 
By direct computation with some substitutions, we see that 
	\bea
		\vek r^{(j)} &=& \eta_{j}{\vek m}_{j+1}^{(j+1)}\nn\\
			 &=& \prn{-s_{j+1}\eta_{j-1}}\prn{-s_{j+1}\vek m_{j}^{(j)} + c_{j+1}\vek v_{j+1}}\nn\\
			 &=&  s_{j+1}^{2}\eta_{j-1}\vek m_{j}^{(j)} + c_{j+1}s_{j+1}\eta_{j-1}\vek v_{j+1}\nn\\
			 &=& \eta_{j-1}\vek m_{j}^{(j)}-c_{j+1}^{2}\eta_{j-1} \vek m_{j}^{(j)} + c_{j+1}s_{j+1}\eta_{j-1}\vek v_{j+1}\nn\\
			 &=& \vek r^{(j-1)}-c_{j+1}\prn{c_{j+1}\eta_{j-1}\vek m_{j}^{(j)} + s_{j+1}\eta_{j-1}\vek v_{j+1}}\nn\\
			 &=& \vek r^{(j-1)} - c_{j+1}\eta_{j-1}\vek m_{j}^{(j+1)}.\nn
	\eea
	Thus the residual subvectors $\vek r_{\vek u}^{(j)}$ and $\vek r_{\vek p}^{(j)}$ satisfy the progressive updating formulas
	\be\nn
	\vek r_{\vek u}^{(j)} = \vek r_{\vek u}^{(j-1)} - c_{j+1}\eta_{j-1}\vek m_{j,\vek u}^{(j+1)}\mand \vek r_{\vek p}^{(j)} = \vek r_{\vek p}^{(j-1)} - c_{j+1}\eta_{j-1}\vek m_{j,\vek p}^{(j+1)}.
	\ee
However, this representation
would require storage of two full-length vectors.

We have seen that through recursion formulas, we can compute one or both subvectors of the residual vector produced by
preconditioned MINRES.  However, we wish to avoid storing or re-evaluating these vectors since we actually want to compute only the appropriate \emph{norm} of each piece.  The full residual minimization is 
with respect to the $\vek P^{-1}$ norm $\norm{\cdot}_{\vek P^{-1}}$.  We have assumed in this paper that the preconditioner has block 
diagonal structure; and as $\vek P$ is symmetric positive definite, it follows that its blocks are as well. Thus we can write
\be\nn
	\bignorm{\vek r^{(j)}}_{\vek P^{-1}}^{2} = \bignorm{\vek r_{\vek u}^{(j)}}_{\vek P_{\vek u}^{-1}}^{2} + \bignorm{\vek r^{(j)}_{\vek p}}_{\vek P_{\vek p}^{-1}}^{2}.
\ee
These norms can also be computed recursively using information already on hand. 
We note that doing it as in the following allows us to avoid additional applications of the preconditioner.

In order to compute the norms separately, it is necessary
that we also use the same fact for the newest Lanczos vector,
\be\nn
	1 = \norm{\vek v_{j+1}}_{\vek P^{-1}}^{2} = \norm{\vek v_{j+1,\vek u}}_{\vek P_{\vek u}^{-1}}^{2} + \norm{\vek v_{j+1,\vek p}}_{\vek P_{\vek p}^{-1}}^{2}.
\ee
Let 
\bea
	\psi_{j+1,\vek u} &=& \norm{\vek v_{j+1,\vek u}}_{\vek P_{\vek u}^{-1}}^{2} = \ip{\vek z_{j+1,\vek u},\vek v_{j+1,\vek u}}\mand\nn\\
	\psi_{j+1,\vek p} &=& \norm{\vek v_{j+1,\vek p}}_{\vek P_{\vek p}^{-1}}^{2} = \ip{\vek z_{j+1,\vek p},\vek v_{j+1,\vek p}}\nn.
\eea
We can store these values so they are available for later use.  We further denote $\eta_{j,\vek u} = \bignorm{\vek r_{\vek u}^{(j)}}_{\vek P_{\vek u}^{-1}}$ and 
$\eta_{j,\vek p} = \bignorm{\vek r_{\vek p}^{(j)}}_{\vek P_{\vek p}^{-1}}$.
This allows us to concisely state the following theorem.
\bthm
	The ratio between the full preconditioned residual norm and the respective preconditioned residual subvector norms 
	can be written progressively as follows:
	\bea
	\prn{\frac{\eta_{j,\vek u}}{\eta_{j}} }^{2} &=& {s_{j+1}^{2}\prn{\frac{\eta_{j-1,\vek u}}{\eta_{j-1}} }^{2} - 2s_{j+1}c_{j+1}\prn{\vek m_{j,\vek u}^{(j)}}^{T}\vek z_{j+1,\vek u} + c_{j+1}^{2}\psi_{j+1,\vek u}}\mand\nn\\
		\prn{\frac{\eta_{j,\vek p}}{\eta_{j}} }^{2} &=& {s_{j+1}^{2}\prn{\frac{\eta_{j-1,\vek p}}{\eta_{j-1}} }^{2} - 2s_{j+1}c_{j+1}\prn{\vek m_{j,\vek p}^{(j)}}^{T}\vek z_{j+1,\vek p} + c_{j+1}^{2}\psi_{j+1,\vek p}}\nn.
	\eea
\ethm
\bproof
The proofs for each subvector are identical, so we prove it only for one.  We can verify  the identity by direct calculation and then by showing
that the quantities arising from the calculation are available by recursion.  We compute
\bea
	\bignorm{\vek r_{\vek u}^{(j)}}_{\vek P_{\vek u}^{-1}}^{2} &=& \bigprn{\eta_{j}\vek m_{j+1,\vek u}^{(j+1)}}^{T}\vek P_{\vek u}^{-1}\bigprn{\eta_{j}\vek m_{j+1,\vek u}^{(j+1)}}\nn\\
				&=& \eta_{j}^{2}\bigprn{-s_{j+1}\vek m_{j,\vek u}^{(j)} + c_{j+1}\vek v_{j+1,\vek u}}^{T}\vek P_{\vek u}^{-1}\bigprn{-s_{j+1}\vek m_{j,\vek u}^{(j)} + c_{j+1}\vek v_{j+1,\vek u}}\nn\\
				&=& \eta_{j}^{2}\Bigprn{s_{j+1}^{2}\bigprn{\vek m_{j,\vek u}^{(j)}}^{T}\vek P_{\vek u}^{-1}\vek m_{j,\vek u}^{(j)} - 2s_{j+1}c_{j+1}\bigprn{\vek m_{j,\vek u}^{(j)}}^{T}\vek P_{\vek u}^{-1}\vek v_{j+1,\vek u} \nn\\
				& & \quad {}+ c_{j+1}^{2}\vek v_{j+1,\vek u}^{T}\vek P_{\vek u}^{-1}\vek v_{j+1,\vek u}}\nn\\
				&=& \eta_{j}^{2}\prn{s_{j+1}^{2}\frac{\eta_{j-1,\vek u}^{2}}{\eta_{j-1}^{2}} - 2s_{j+1}c_{j+1}\prn{\vek m_{j,\vek u}^{(j)}}^{T}\vek z_{j+1,\vek u} + c_{j+1}^{2}\psi_{j+1,\vek u}}\nn,
\eea
\eproof
Basic algebra finishes the proof.

For completeness, we now present a modified version of \cite[Algorithm~4.1]{ESW-Book.2014}, here called
\cref{alg.prec.minres-modified}.  Note that we omit superscripts for 
$\vek m_{i,\vek p}^{(\ell)}$ and some other subscripts where they are not needed.  Furthermore, we introduce the
scalars $\theta_{\vek u}$ and $\theta_{\vek p}$ defined by 
\be\label{eq:definition_of_theta}
\theta_{\vek u} = \bigip{\vek m_{j,\vek u}^{(j)},\vek z_{j+1,\vek u}}\mand\theta_{\vek p} = \bigip{\vek m_{j,\vek p}^{(j)}, \vek z_{j+1,\vek p}}
\ee
and the squared residual norm fractions
\be\label{eq:definition_of_mu}
\mu_{\vek u} = \prn{\frac{\eta_{j,\vek u}}{\eta_{j}} }^{2}\mand
\mu_{\vek p} = \prn{\frac{\eta_{j,\vek p}}{\eta_{j}} }^{2}.
\ee
Notice that one can work with full-length vectors and the partitioning is only important for the computation of partial inner products as in \eqref{eq:definition_of_theta}.

\balg
\caption{The preconditioned MINRES method with monitoring of $\abs{\eta_{j,\vek u}} = \bignorm{\vek r_{\vek u}^{(j)}}_{\vek P_{\vek u}^{-1}}$ and $\abs{\eta_{j,\vek p}} = \bignorm{\vek r_{\vek p}^{(j)}}_{\vek P_{\vek p}^{-1}}$}\label{alg.prec.minres-modified}
\SetKwInOut{Input}{Input}\SetKwInOut{Output}{Output}
\Input{$\vek K\in \R^{(m+p)\times (m+p)}$, symmetric, and $\vek f \in \R^{m+p}$ as in \eqref{eqn.ABupf}, $\vek P = \blkdiag(\vek P_{\vek u},\vek P_{\vek p})\in\R^{(m+p)\times (m+p)}$, symmetric positive definite.}
\Output{$\vek x^{(j)}$ for some $j$ such that $\norm{\vek r^{(j)}}_{\vek P^{-1}}$ satisfies some convergence criteria.}
	Set $\vek v_{0} = \vek 0$, \quad $\vek w_{0} = \vek w_{1} = \vek 0$\\
	Choose $\vek x^{(0)}$, \quad set $\vek v_{1} = \bbmat \vek f_{1}\\ \vek f_{2} \ebmat - \vek K \, \vek x^{(0)}$, \quad $\vek z_{1} = \vek P^{-1}\vek v_{1}$\\ 
	$\gamma_{1} = \sqrt{\ip{\vek z_{1},\vek v_{1}}}$, \quad $\vek v_{1}\leftarrow \vek v_{1}/\gamma_{1}$, \quad $\vek z_{1}\leftarrow \vek z_{1}/\gamma_{1}$\\ 
	$\psi_{\vek u}=\ip{\vek z_{1,\vek u},\vek v_{1,\vek u}}$, \quad $\psi_{\vek p}=\ip{\vek z_{1,\vek p},\vek v_{1,\vek p}}$\\ 
	$\mu_{\vek u} = \psi_{\vek u}$, \quad $\mu_{\vek p} = \psi_{\vek p}$\\
	$\vek m_{1} = \vek v_{1}$\\
	Set $\eta_{0} = \gamma_{1}$, \quad $\eta_{0,\vek u} = \gamma_{1}\, \sqrt{\psi_{\vek u}}$, \quad $\eta_{0,\vek p} = \gamma_{1}\, \sqrt{\psi_{\vek p}}$, \quad  $s_{0} = s_{1} = 0$, $c_{0}=c_{1} = 1$\\
	\For{$j=1$ until convergence}{
		$\delta_{j} = \ip{\vek K\vek z_{j},\vek z_{j}}$\\
		$\vek v_{j+1} = \vek K\vek z_{j} - \delta_{j}\vek v_{j} - \gamma_{j}\vek v_{j-1}$ \tcp*{Lanczos}
		$\vek z_{j+1} = \vek P^{-1}\vek v_{j+1}$\\
		$\gamma_{j+1} = \sqrt{\ip{\vek z_{j+1},\vek v_{j+1}}}$\\
		$\vek v_{j+1} \leftarrow \vek  v_{j+1}/\gamma_{j+1}$\\
		$\vek z_{j+1} \leftarrow \vek  z_{j+1}/\gamma_{j+1}$\\
		$\alpha_{0} = c_{j}\delta_{j} - c_{j-1}s_{j}\gamma_{j}$ \tcp*{Update QR factorization}
		$\alpha_{1} = \sqrt{\alpha_{0}^{2} + \gamma_{j+1}^{2}}$\\
		$\alpha_{2} = s_{j}\delta_{j} + c_{j-1}c_{j}\gamma_{j}$\\
		$\alpha_{3} = s_{j-1}\gamma_{j}$\\
		$c_{j+1} = \alpha_{0}/\alpha_{1}$, \quad $s_{j+1}=\gamma_{j+1}/\alpha_{1}$ \tcp*{Givens rotations}
		$\theta_{\vek u} \leftarrow \ip{\vek m_{j,\vek u}, \vek z_{j+1,\vek u}}$, \quad $\theta_{\vek p} \leftarrow \ip{\vek m_{j,\vek p}, \vek z_{j+1,\vek p}}$\\
		$\psi_{\vek u} \leftarrow \ip{\vek z_{j+1,\vek u},\vek v_{j+1,\vek u}}$, \quad $\psi_{\vek p} \leftarrow \ip{\vek z_{j+1,\vek p},\vek v_{j+1,\vek p}}$\\
		$\vek m_{j+1} = -s_{j+1}{\vek m}_{j} + c_{j+1}\vek v_{j+1}$\\
		$\vek w_{j+1} = (\vek z_{j} - \alpha_{3}\vek w_{j-1} - \alpha_{2}\vek w_{j})/\alpha_{1}$\\
		$\vek x^{(j)} = \vek x^{(j-1)} + c_{j+1}\eta_{j-1}\vek w_{j+1}$\\
		$\mu_{\vek u} \leftarrow s_{j+1}^{2}\mu_{\vek u} - 2s_{j+1}c_{j+1}\theta_{\vek u} + c_{j+1}^{2}\psi_{\vek u}$\\
		$\mu_{\vek p} \leftarrow s_{j+1}^{2}\mu_{\vek p} - 2s_{j+1}c_{j+1}\theta_{\vek p} + c_{j+1}^{2}\psi_{\vek p}$\\
		$\eta_{j} = -s_{j+1}\eta_{j-1}$\tcp*{total residual norm}
		$\eta_{j,\vek u} = \eta_{j}\, \sqrt{\mu_{\vek u}}$, \quad $\eta_{j,\vek p} = \eta_{j}\, \sqrt{\mu_{\vek p}}$\tcp*{partial residual norms}
		Test for convergence
	}
\ealg

\section{Examples and Numerical Experiments}\label{section.examples}

In the introduction, we motivated our study by mentioning that the residual subvectors in saddle-point systems \eqref{eqn.ABupf} often have different physical interpretations.
In this section, we discuss several examples to underline this fact, and we also present numerical results which demonstrate the correctness of \cref{alg.prec.minres-modified}.
We do this by storing all the residual vectors throughout the iteration and evaluating the preconditioned subvector norms a-posteriori.
This is a debugging step which is introduced to verify the correctness of \cref{alg.prec.minres-modified}.
Our \matlab\ implementation of \cref{alg.prec.minres-modified} is available at \cite{SUBMINRES}.

In all examples, we begin with an all-zero initial guess and we stop when the relative reduction of the total residual 
\begin{equation*}
	\frac{\eta_j}{\eta_0} 
	= \frac{\bignorm{\vek r^{(j)}}_{\vek P^{-1}}}{\bignorm{\vek r^{(0)}}_{\vek P^{-1}}}
\end{equation*}
falls below $10^{-6}$.
We could easily utilize a refined stopping criterion such as
\begin{equation*}
	\abs{\eta_{j,\vek u}} = \bignorm{\vek r_{\vek u}^{(j)}}_{\vek P_{\vek u}^{-1}} \le \varepsilon_{\vek u}
	\quad \text{and} \quad 
	\abs{\eta_{j,\vek p}} = \bignorm{\vek r_{\vek p}^{(j)}}_{\vek P_{\vek p}^{-1}} \le \varepsilon_{\vek p}
\end{equation*}
to take advantage of the ability to monitor the residual subvector norms.
Concrete applications of this (e.g., in optimization) are outside the scope of this paper and will be addressed elsewhere.

\begin{figure}[htbp]
	\centering
	\includegraphics[width=0.45\textwidth]{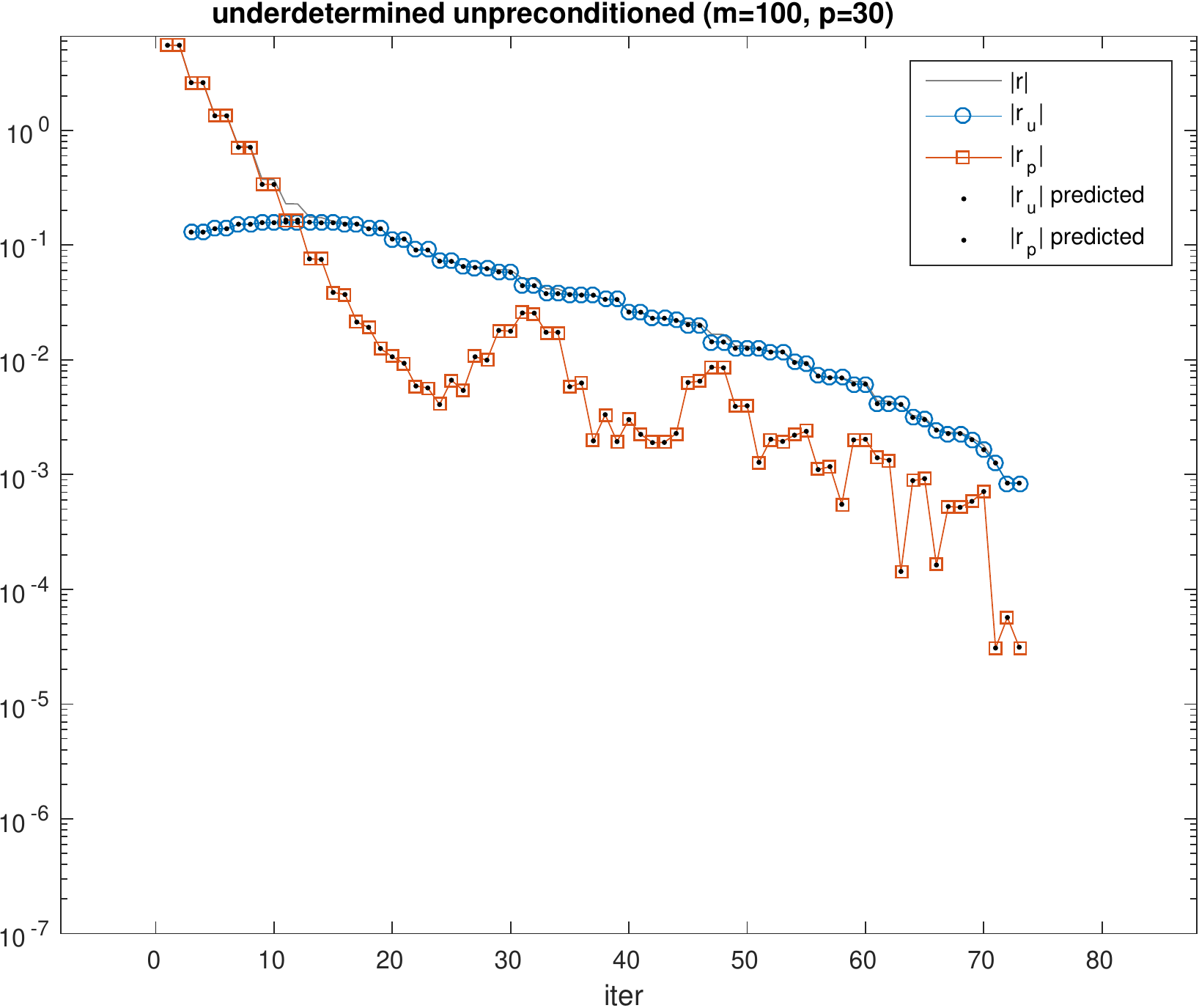}
	\includegraphics[width=0.45\textwidth]{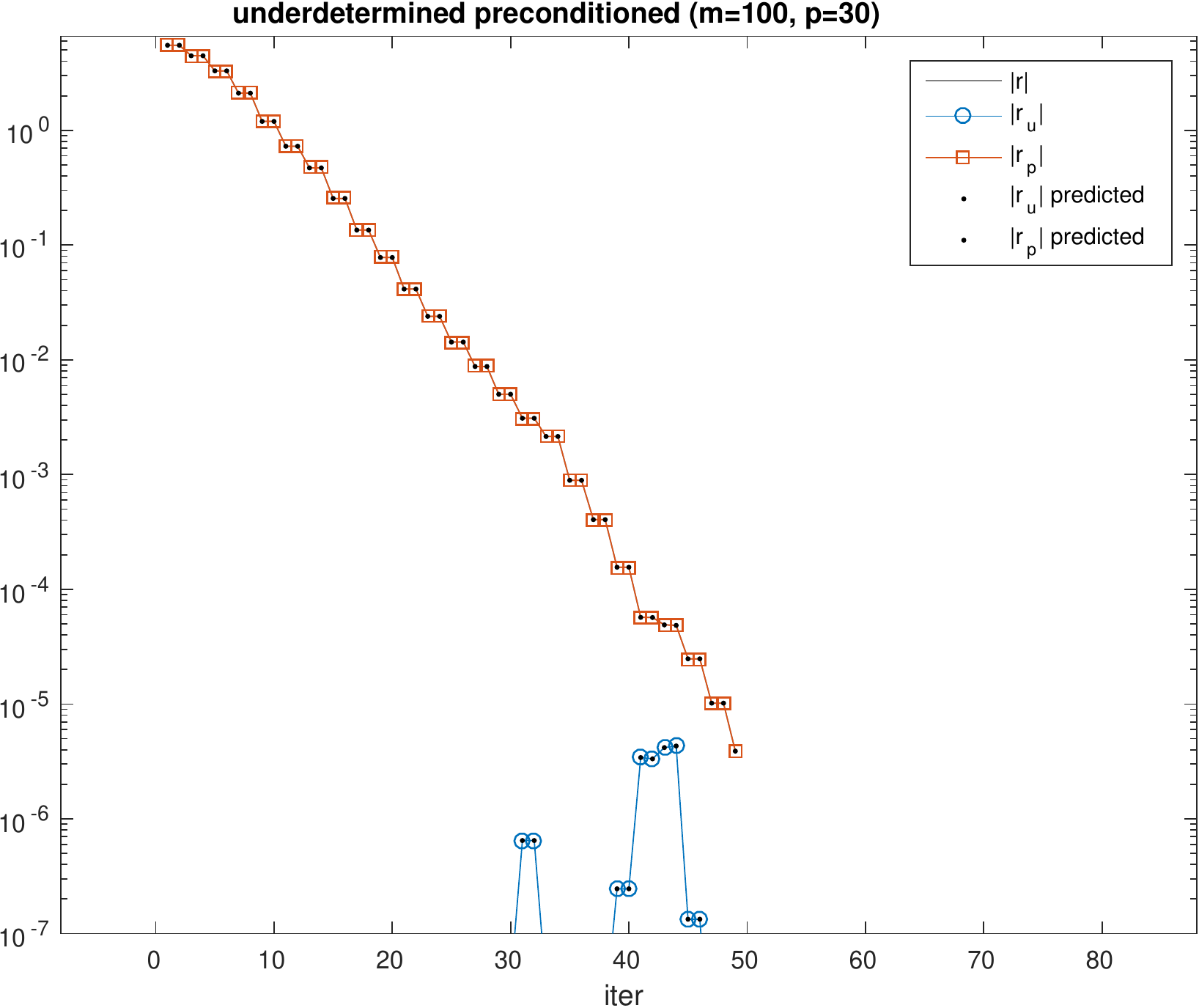}
	\caption{Convergence history of the residual subvectors for \cref{example.underdetermined} (underdetermined linear system) in the unpreconditioned (left) and preconditioned case (right). The plot, as all following convergence plots, also confirms the correctness of \cref{alg.prec.minres-modified} and of our implementation. Notice that in the right figure, the first residual norm $\norm{\vek r_{\vek u}}_{\vek P_{\vek u}^{-1}}$ is partially outside of the plot range. This results from our choice to maintain the same scales for both figures.}
	\label{fig:underdetermined}
\end{figure}

\begin{example}[Least-Norm Solution of Underdetermined Linear System] \hfill \\
	\label{example.underdetermined}
	We consider an underdetermined and consistent linear system $\vek B \vek u = \vek b$ with a matrix $\vek B \in \R^{m \times n}$ of full row rank and $m < n$.
	Its least-norm solution 
	\begin{equation*}
		\text{Minimize} \quad \frac{1}{2} \norm{\vek u}_{\vek H}^2 \msuchthat \vek B \vek u = \vek b, \quad \vek u \in \R^n
	\end{equation*}
	in the sense of the inner product defined by the symmetric positive definite matrix $\vek H \in \R^{n \times n}$, is uniquely determined by the saddle-point system
	\begin{equation*}
		\bbmat \vek H & \vek B^{T}\\ \vek B & \Bzero \ebmat \bbmat \vek u \\ \vek p \ebmat = \bbmat \vek \Bzero\\ \vek b \ebmat.
	\end{equation*}
	Notice that the first block of equations represents optimality, while the second block represents feasibility of a candidate solution $(\vek u, \vek p)^{T}$. \\
	Our test case uses data created through the commands
	\begin{lstlisting}[basicstyle={\ttfamily\upshape},gobble=6]
	rng(42); n=100; m=30; B=randn(m,n); b=randn(m,1);
	H=spdiags(rand(n,1),0,n,n); 
	\end{lstlisting}
	As preconditioner we use either $\vek P_1 = \blkdiag(\vek I_{n \times n}, \vek I_{m \times m})$ (the unpreconditioned case) or $\vek P_2 = \blkdiag(\vek H, \vek I_{m \times m})$.
	The convergence histories in \cref{fig:underdetermined} show that the amount by which the two residual subvectors contribute to their combined norm may indeed be quite different, and it depends on the preconditioner. 
	In this example, the feasibility residual $\vek r_{\vek p} = \vek b - \vek B \vek u$ is significantly smaller than the optimality residual $\vek r_{\vek u} = - \vek H \vek u - \vek B^{T} \vek p$ in the unpreconditioned case where we used $\vek P_1$.
	To be more precise, the average value of $\mu_{\vek p}$ throughout the iteration history is about 21\%.
	Quite the opposite is true for the case of $\vek P_2$, when $\mu_{\vek p}$ is close to 100\%.
\end{example}

\begin{figure}[htbp]
	\centering
	\includegraphics[width=0.45\textwidth]{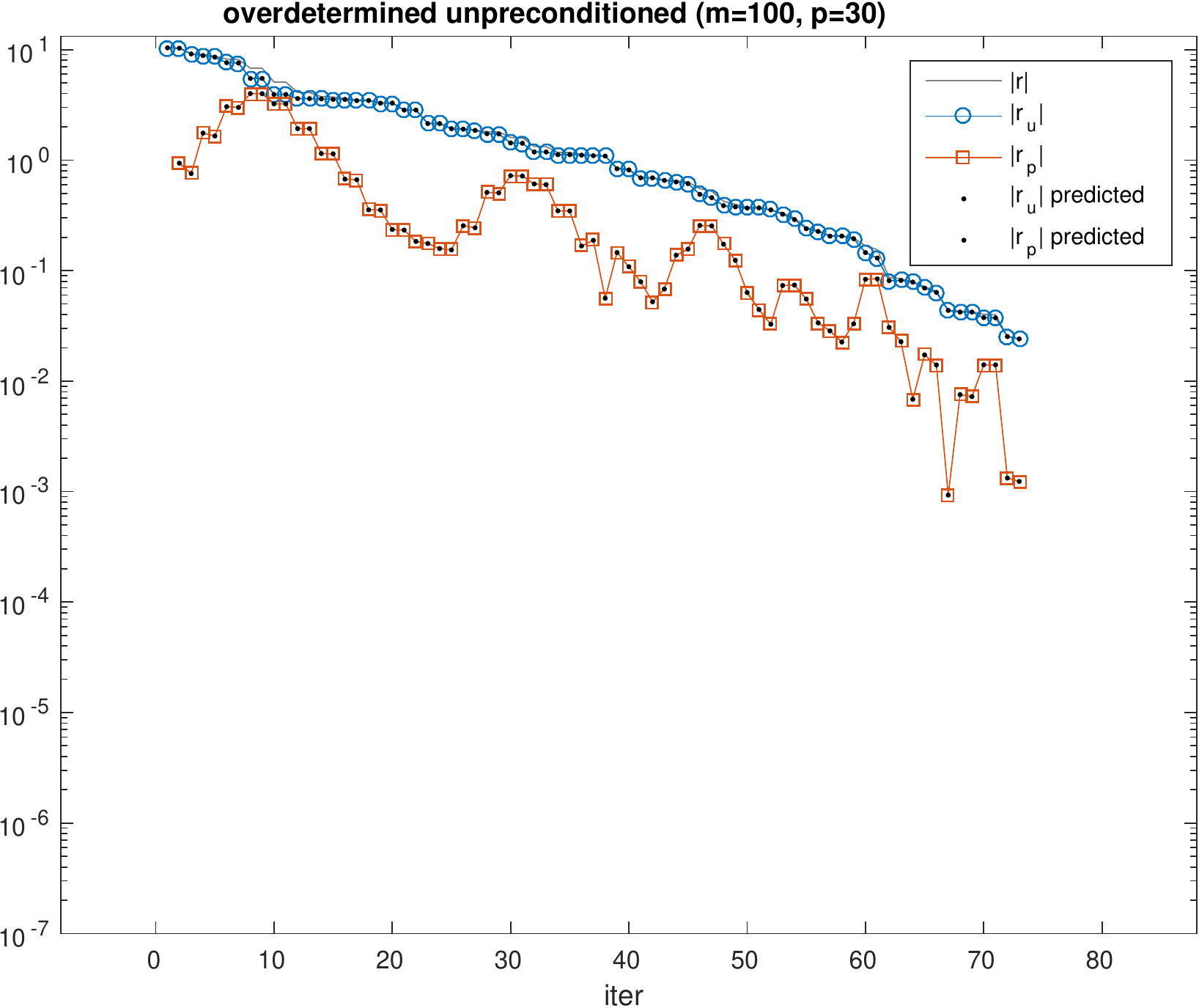}
	\includegraphics[width=0.45\textwidth]{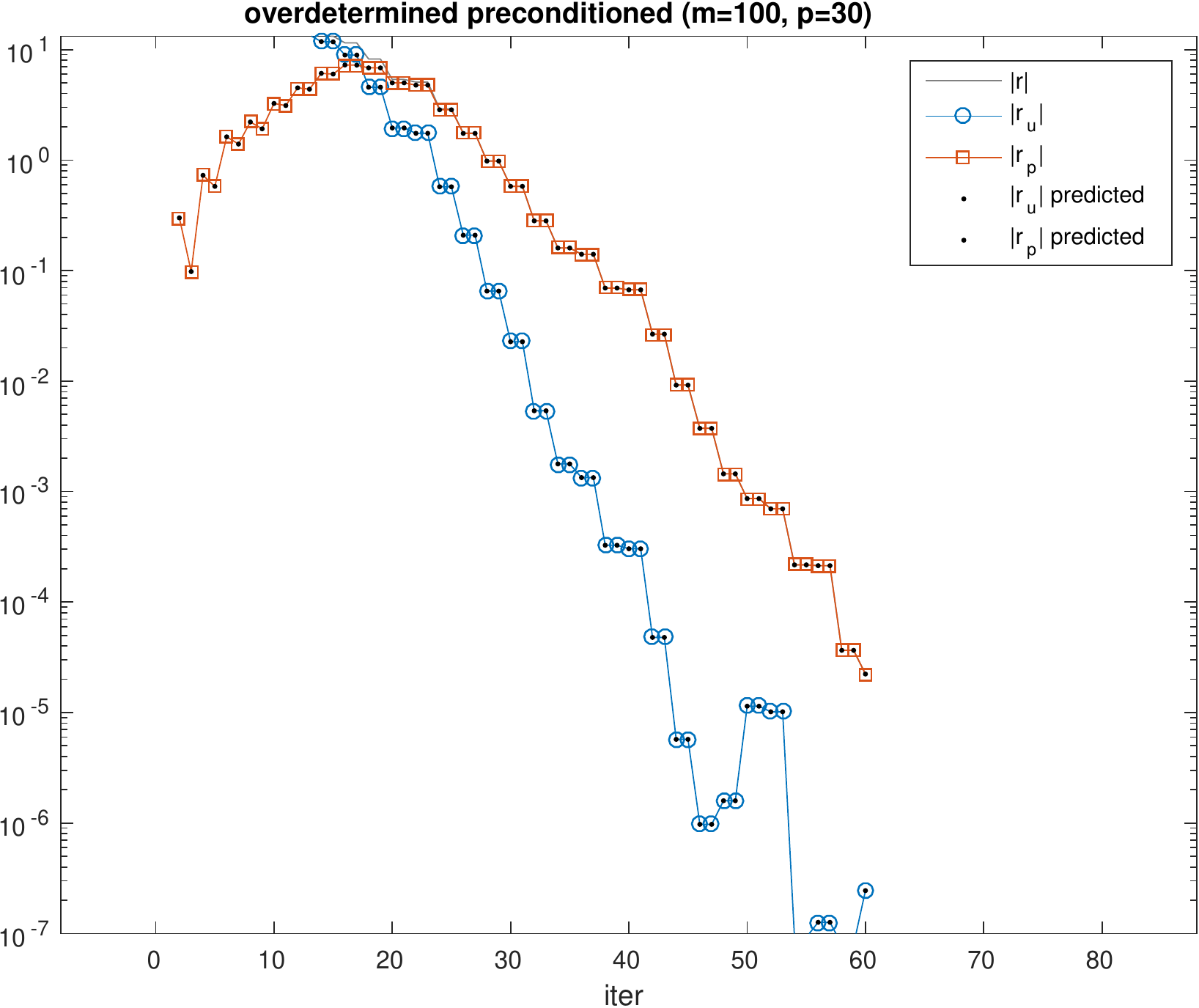}
	\caption{Convergence history of the residual subvectors for \cref{example.overdetermined} (overdetermined linear system) in the unpreconditioned (left) and preconditioned case (right).}
	\label{fig:overdetermined}
\end{figure}

\begin{example}[Least-Squares Solution of Overdetermined Linear System] \hfill \\
	\label{example.overdetermined}
	Here we consider an overdetermined linear system $\vek B^{T} \vek p = \vek b$ with a matrix $\vek B \in \R^{m \times n}$ of full row rank and $m < n$.
	Its least-squares solution 
	\begin{equation*}
		\text{Minimize} \quad \frac{1}{2} \norm{\vek B^{T} \vek p - \vek b}_{\vek H^{-1}}^2, \quad \vek p \in \R^m
	\end{equation*}
	with $\vek H \in \R^{n \times n}$ symmetric positive definite as above, is uniquely determined by the normal equations, $\vek B \, \vek H^{-1} (\vek B^{T} \vek p - \vek b) = \Bzero$.
	By defining $\vek u$ as the 'preconditioned residual' $\vek H^{-1} (\vek B^{T} \vek p - \vek b)$, we find that the least-squares solution is in turn equivalent to the saddle-point system
	\begin{equation*}
		\bbmat \vek H & \vek B^{T}\\ \vek B & \Bzero \ebmat \bbmat \vek u \\ \vek p \ebmat = \bbmat \vek b\\ \vek \Bzero \ebmat.
	\end{equation*}
	The first block of equations now represents feasibility for the constraint defining the auxiliary quantity $\vek u$, while the second requires $\vek u$ to lie in the kernel of $B$. \\
	Similarly as above, we derive test data through
	\begin{lstlisting}[basicstyle={\ttfamily\upshape},gobble=6]
	rng(42); n=100; m=30; B=randn(m,n); b=randn(n,1);
	H=spdiags(rand(n,1),0,n,n);
	\end{lstlisting}
	We employ the same preconditioners $\vek P_1$ and $\vek P_2$ as in \cref{example.underdetermined}.
	Once again, the convergence behavior of the two residual subvectors is fundamentally different (\cref{fig:overdetermined}) for the two cases: in the unpreconditioned case, the average value of $\mu_{\vek p}$ is about 9\%, while it is 72\% in the preconditioned case.
\end{example}

Our remaining examples involve partial differential equations and they will be stated in variational form, by specifying bilinear forms $a(\cdot,\cdot): V \times V \to \R$, $b(\cdot,\cdot): V \times Q \to \R$, and, where appropriate, $c(\cdot,\cdot): Q \times Q \to \R$.
Here $V$ and $Q$ are (real) Hilbert spaces.
The matrices $\vek A$, $\vek B$ and $\vek C$ in \eqref{eqn.ABupf} are then obtained by evaluating the respective bilinear forms on a basis of an appropriate finite dimensional subspace $V_h$ or $Q_h$, e.g., $[\vek B]_{ij} = b(\varphi_j,\psi_i)$ for basis elements $\varphi_j \in V_h$ and $\psi_i \in Q_h$.
Similarly, the right hand side vectors $\vek f_{\vek u}$ and $\vek f_{\vek p}$ are obtained from evaluating problem dependent linear forms $f_{\vek u}(\cdot)$ and $f_{\vek p}(\cdot)$ on the same basis functions $\varphi_i$ and $\psi_i$, respectively.
Finally, the matrices and vectors obtained in this way may need to be updated due to the incorporation of essential (Dirichlet) boundary conditions.

\begin{figure}[htbp]
	\centering
	\includegraphics[width=0.45\textwidth]{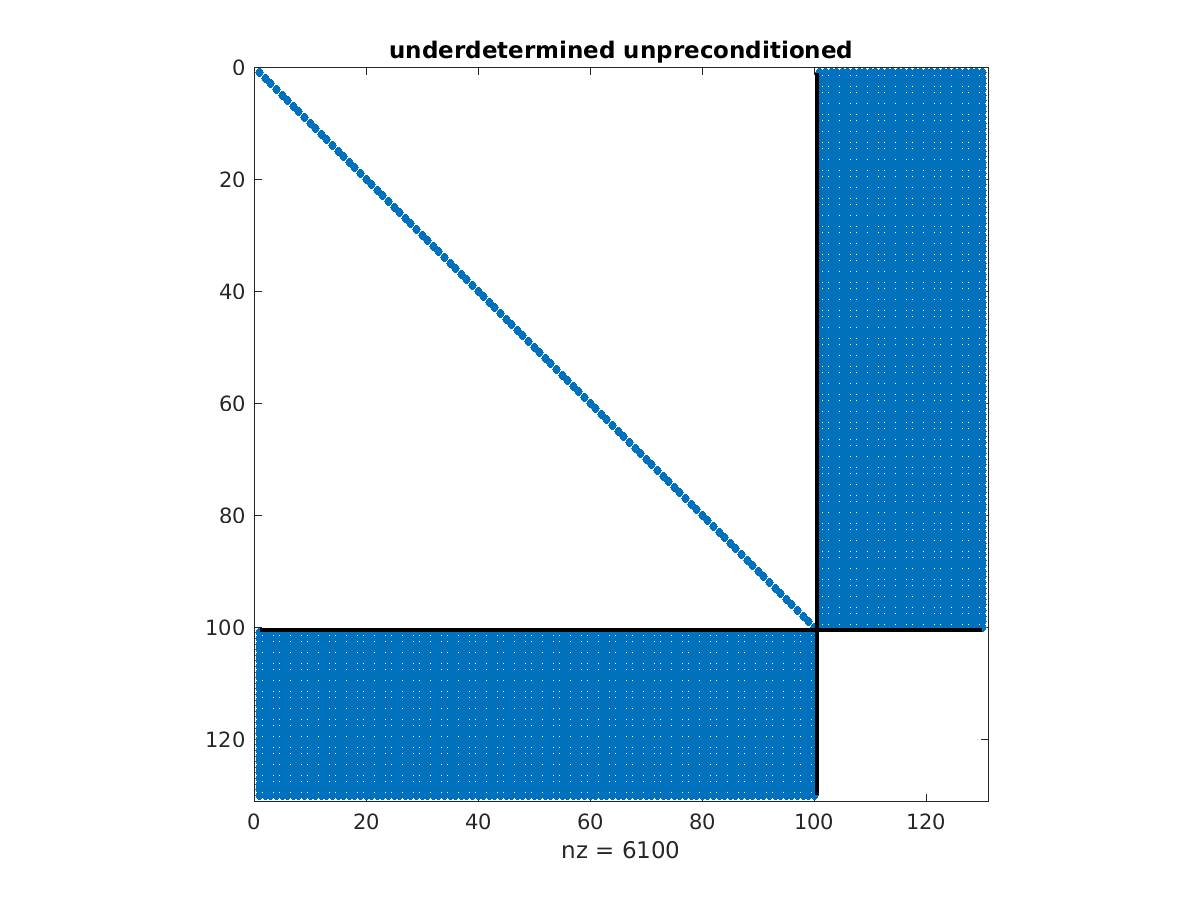}
	\includegraphics[width=0.45\textwidth]{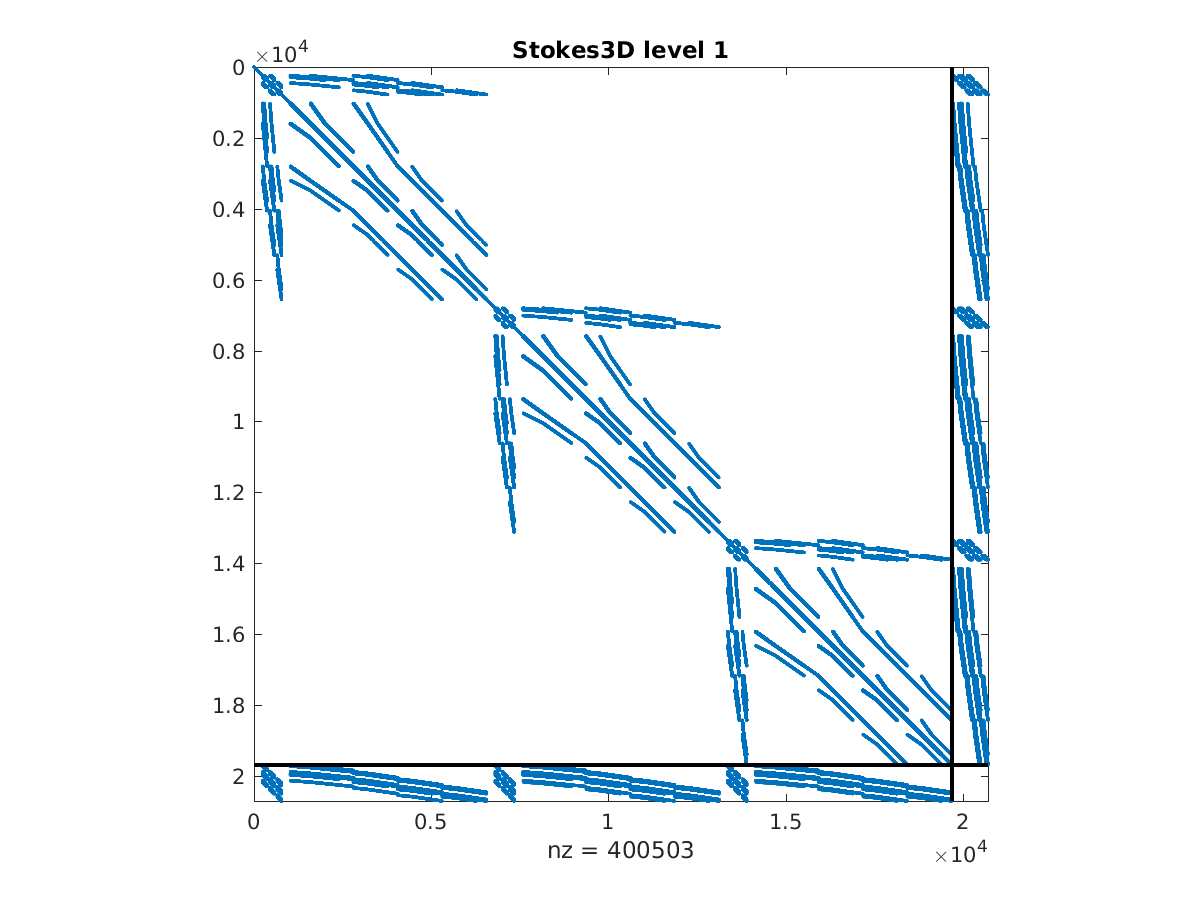} \\
	\includegraphics[width=0.45\textwidth]{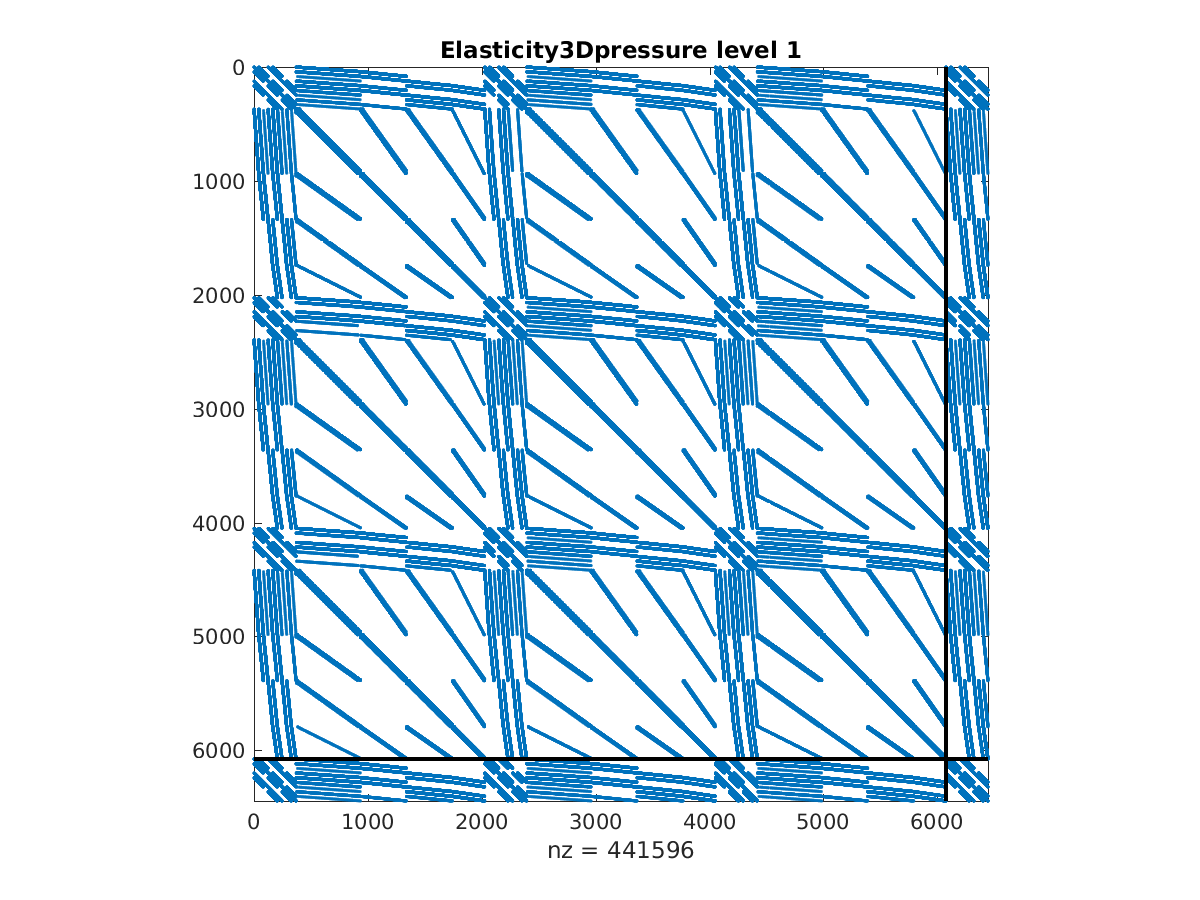}
	\includegraphics[width=0.45\textwidth]{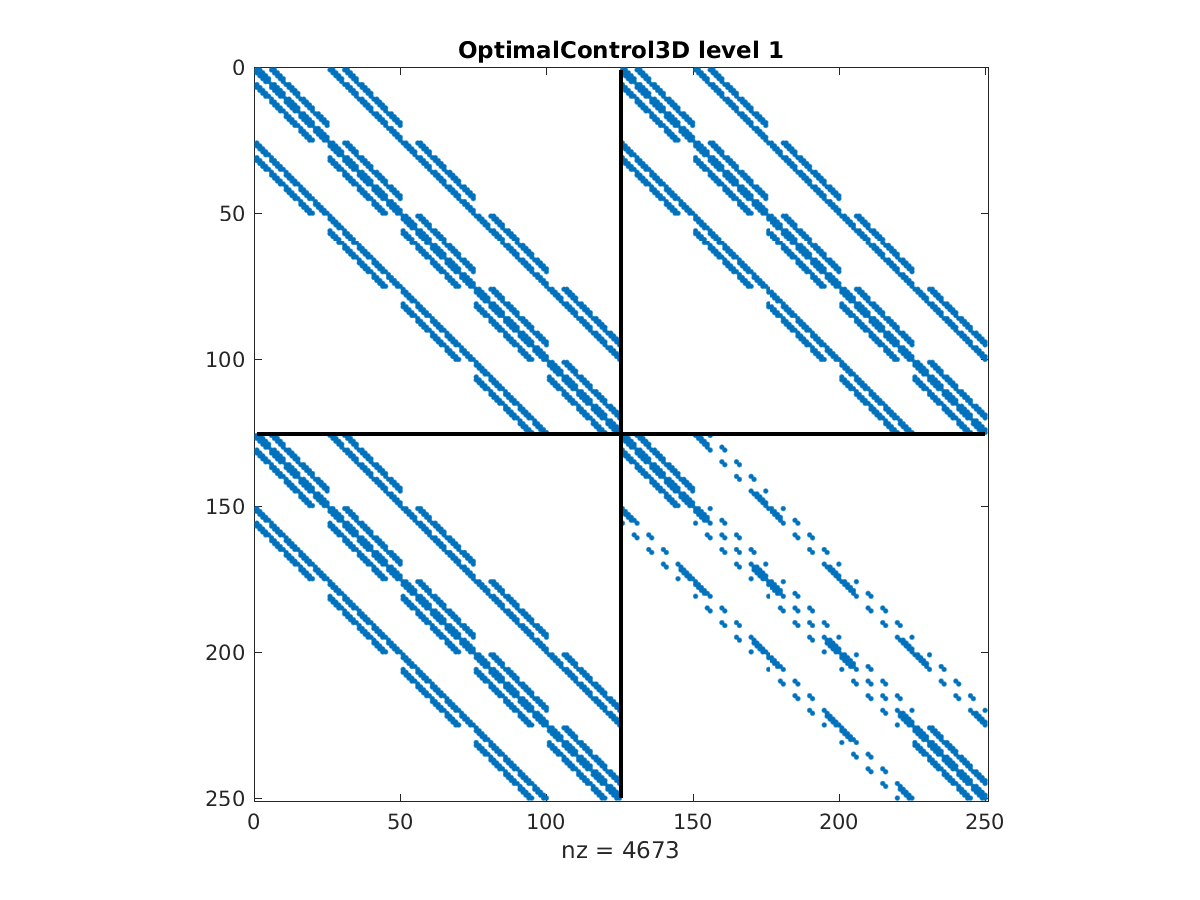}
	\caption{Sparsity plots of the saddle-point systems arising in \cref{example.underdetermined} (underdetermined linear system), \cref{example.Stokes} (Stokes), \cref{example.Elasticity} (elasticity), and \cref{example.OptimalControl} (optimal boundary control).}
	\label{fig:spy}
\end{figure}

All numerical tests were conducted using the \python\ interface of the finite element library \fenics\ \cite{LoggMardalWells2012:1} (version~1.6) to generate the matrices and vectors.
Those were then exported in \petsc\ binary format and read from \matlab\ through the helper functions \file{readPetscBinMat.m} and \file{readPetscBinVec.m} provided by Samar Khatiwala on his web page.
We deliberately turned off the reordering feature of \fenics\ for the degrees of freedom to preserve the block structure of \eqref{eqn.ABupf} for illustration purposes (see \cref{fig:spy}) and in order for the subvectors $\vek u$ and $\vek p$ and the residuals $\vek r_{\vek u}$ and $\vek r_{\vek p}$ to remain contiguous in memory.
However, our theory does not rely on a particular ordering of the subvector components, and our implementation of \cref{alg.prec.minres-modified} allows for arbitrary component ordering.

It will turn out to be useful for the following examples to specify the physical units (Newton: \si{\newton}, meters: \si{\meter}, seconds: \si{\second}, Watt: \si{\watt}, and Kelvin: \si{\kelvin}) for all involved quantities in the following examples.

\begin{figure}[htbp]
	\centering
	\includegraphics[width=0.45\textwidth]{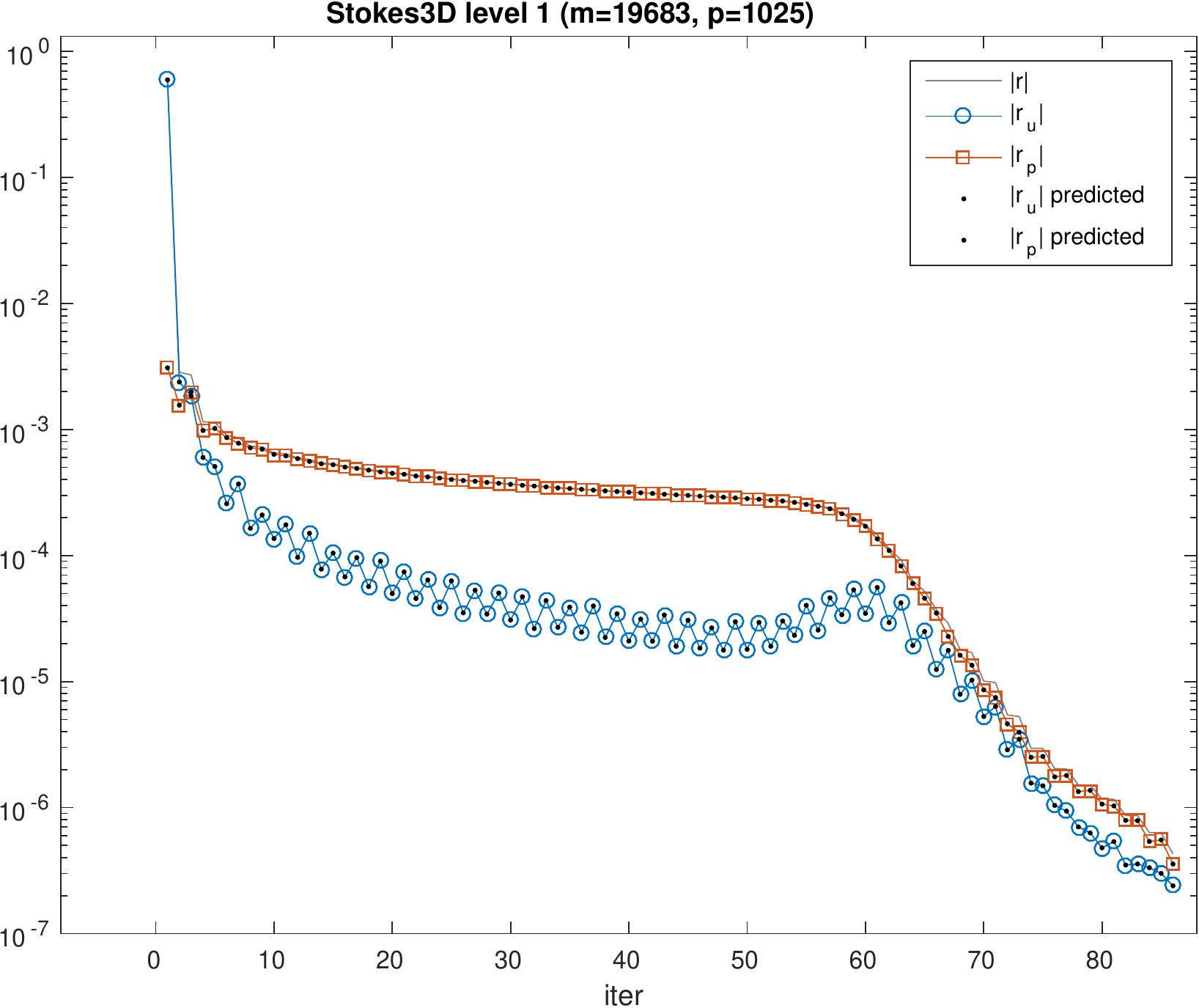}
	\includegraphics[width=0.45\textwidth]{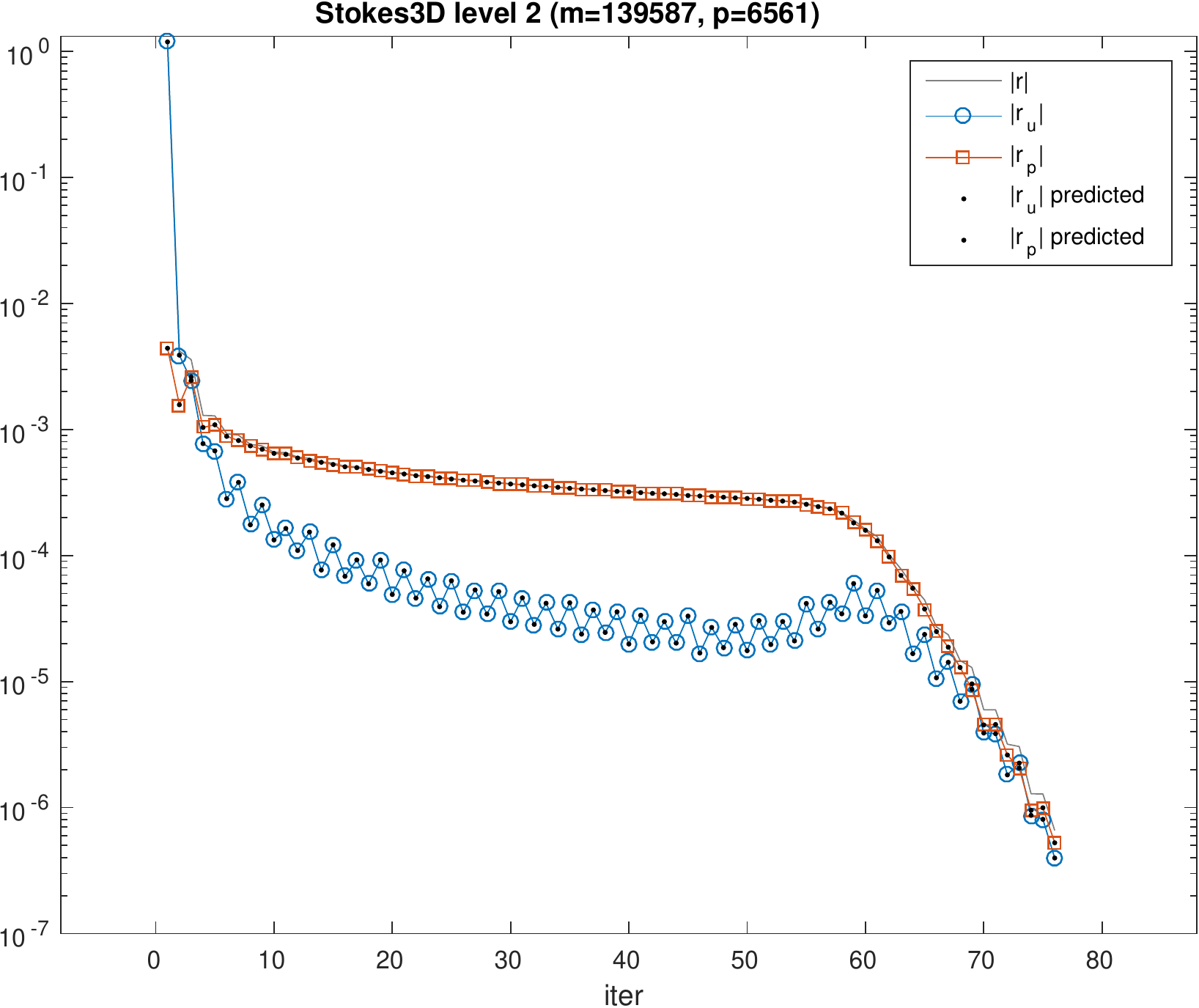}
	\caption{Convergence history of the residual subvectors for \cref{example.Stokes} (Stokes) on a coarse grid (left) and its uniform refinement (right).}
	\label{fig:Stokes}
\end{figure}

\begin{example}[Stokes Channel Flow] \hfill \\
	\label{example.Stokes}
	We consider the variational formulation of a 3D stationary Stokes channel flow configuration within the domain $\Omega = (0,10) \times (0,1) \times (0,1)$.
	Dirichlet conditions are imposed on the fluid velocity $\vek u$ everywhere except at the 'right' (outflow) boundary $(x = 10)$, where do-nothing conditions hold.
	The Dirichlet conditions are homogeneous (no-slip) except at the 'left' (inflow) boundary $(x = 0)$, where $\vek u(x,y,z) = \vek u_\textup{in}(x,y,z) = (y \, (1-y) \, z \, (1-z),0,0)^{T}$~\si{\meter\per\second} is imposed.

	Appropriate function spaces for this setup are $V = \{ \vek v \in H^1(\Omega;\R^3): \vek v = \Bzero \text{ on } \Gamma \setminus \Gamma_\textup{right} \}$ for the velocity and $Q = L^2(\Omega)$ for the pressure.
	The relevant bilinear and linear forms associated with this problem are
	\begin{equation*}
		a(\vek u, \vek v) = \mu \int_\Omega \nabla \vek u \dprod \nabla \vek v \, \dx, \quad
		b(\vek u, q) = \int_\Omega q \, \operatorname{div} \vek u \, \dx \mand 
		f_{\vek u}(\vek v) = \int_\Omega \vek f \cdot \vek v \, \dx.
	\end{equation*}
	We use the dynamic viscosity parameter $\mu = \SI{1e-3}{\newton\second\per\meter\squared}$ (water) and 
	zero
	right hand side force 
	$\vek f = (0,0,0)^{T}$~\si{\newton\per\meter\cubed}.

	By considering units, or by investigating the underlying physics, we infer that the first component of the residual, $r_{\vek u} = f_{\vek u} - a(\vek u,\cdot) - b(\cdot,p)$, represents a net sum of forces, measured in \si{\newton}.
	Similarly, the second residual $r_{\vek p} = - b(\vek u,\cdot)$ represents the net flux of fluid through the impermeable channel walls, measured in \si{\meter\cubed\per\second}.
	Clearly, both parts of the residual must be zero at the converged solution, but their departure from zero at intermediate iterates has different physical interpretations.

	As preconditioner $\vek P = \blkdiag(\vek P_{\vek u},\vek P_{\vek p})$, we use the block diagonal matrix induced by the bilinear forms
	\begin{equation*} 
		a(\vek u,\vek v) \mand \mu^{-1} \int_\Omega p \, q \, \dx,
	\end{equation*}
	respectively, similar to \cite[Section~4.2]{ESW-Book.2014}, where the constant-free problem is considered.
	Notice that the inclusion of the constant $\mu^{-1}$ into the pressure mass matrix renders the preconditioner compatible with the physical units of the problem.

	For our numerical test, we discretized the problem using the Taylor-Hood finite element. 
	The homogeneous and non-homogeneous Dirichlet boundary conditions were included by modifying the saddle-point components $\vek A$, $\vek B$, $\vek B^{T}$, and right hand side $\vek f_{\vek u}$ in \eqref{eqn.ABupf} in a symmetric way through the \lstinline[basicstyle={\ttfamily\upshape}]!assemble_system! call in \fenics.
	Notice that this effectively modifies some components of the first residual subvector $\vek r_{\vek u}$ by expressions of the form $\vek u_\textup{in}(x,y,z)- \vek u(x,y,z)$, measured in \si{\meter\per\second}.
	The same modifications apply to the preconditioner $\vek P$.

	\cref{fig:Stokes} displays the convergence behavior of the residual subvector norms on a relatively coarse mesh and its refinement.
	The result illustrates the mesh independence of the preconditioned iteration, and it also shows that the residual $\vek r_{\vek p}$ (representing mass conservation) lags behind the residual $\vek r_{\vek u}$ over the majority of the iterations.
\end{example}

\begin{figure}[htbp]
	\centering
	\includegraphics[width=0.45\textwidth]{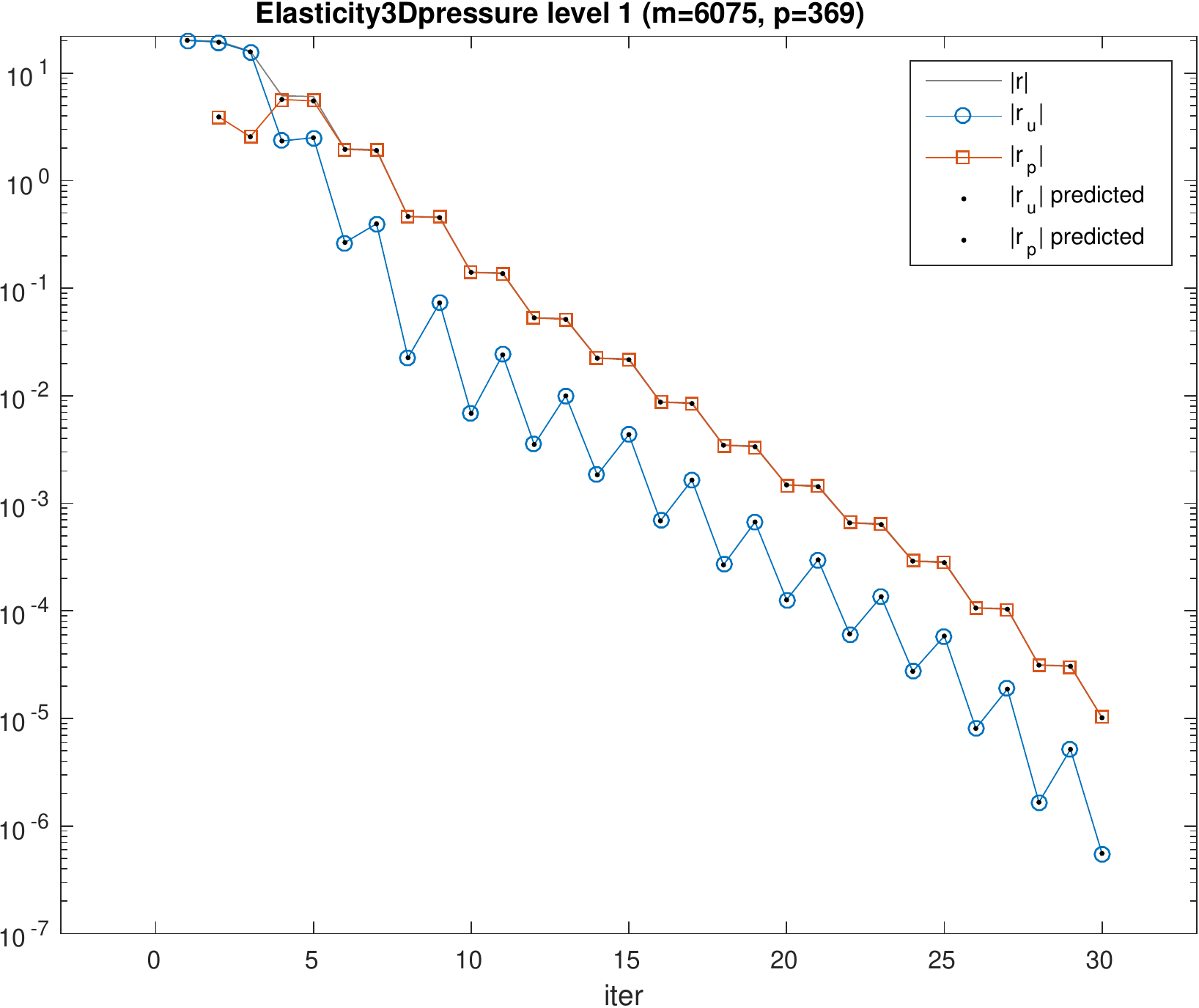}
	\includegraphics[width=0.45\textwidth]{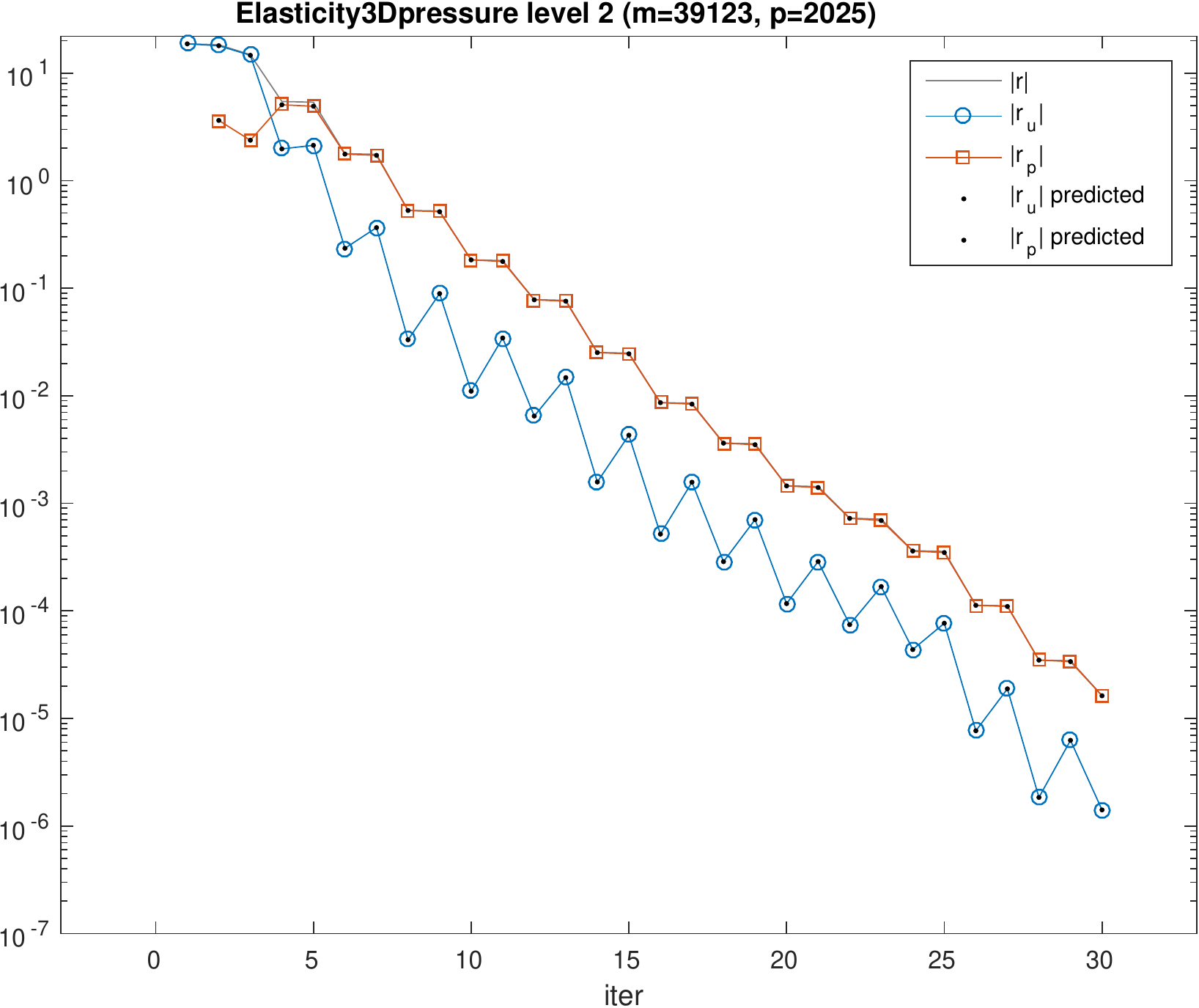}
	\caption{Convergence history of the residual subvectors for \cref{example.Elasticity} (elasticity) on a coarse grid (left) and its uniform refinement (right).}
	\label{fig:Elasticity}
\end{figure}

\begin{example}[Linear Elasticity with Nearly Incompressible Material]  \hfill \\
	\label{example.Elasticity}
	This example describes a tensile test with a rod of square cross section, which occupies the domain $\Omega = (0,100) \times (0,10) \times (0,10)$.
	We use \si{\milli\meter} here in place of \si{\meter} as our length unit.
	Homogeneous Dirichlet conditions for the displacement $\vek u$ are imposed at the 'left' (clamping) boundary $(x = 0)$, while natural (traction) boundary conditions are imposed elsewhere.
	The imposed traction pressure is zero except at the 'right' (forcing) boundary $(x = 100)$, where a uniform pressure of $\vek g = (1,0,0)^{T}$~\si{\newton\per\milli\meter\squared} is imposed.
	As is customary for nearly incompressible material, we introduce an extra variable $p$ for the hydrostatic pressure (see for instance \cite{Wieners2000}) in order to overcome the ill-conditioning of a purely displacement based formulation known as \emph{locking} \cite{BabuskaSuri1992:1}.
	We employ the standard isotropic stress-strain relation, $\boldsymbol{\sigma} = 2 \, \mu \, \boldsymbol{\varepsilon}(\vek u) + p \, \vek I$ and $\operatorname{div} \vek u = \lambda^{-1} p$.
	Here $\boldsymbol{\varepsilon}(\vek u) = (\nabla \vek u + \nabla \vek u^{T})/2$ denotes the symmetrized Jacobian of $\vek u$, while $\mu$ and $\lambda$ denote the Lam{\'e} constants.
	We choose as material parameters Young's modulus $E = \SI{50}{\newton\per\milli\meter\squared}$ and a Poisson ratio of $\nu = 0.49$.
	These particular values describe a material like nearly incompressible rubber, and a conversion to the Lam{\'e} constants yields $\mu = \frac{E}{2 \, (1+\nu)} = \SI{16.78}{\newton\per\milli\meter\squared}$ and $\lambda = \frac{\nu E}{(1+\nu)(1-2\,\nu)} = \SI{822.15}{\newton\per\milli\meter\squared}$.

	The variational mixed formulation obtained in this way is described by the spaces $V = \{ \vek v \in H^1(\Omega;\R^3): \vek v = \Bzero \text{ on } \Gamma_\textup{left} \}$ for the displacement and $Q = L^2(\Omega)$ for the hydrostatic pressure.
	The bilinear and linear forms associated with this problem are
	\begin{equation*}
		\begin{aligned}
			& a(\vek u, \vek v) = 2 \, \mu \int_\Omega \boldsymbol{\varepsilon}(\vek u) \dprod \boldsymbol{\varepsilon} (\vek v) \, \dx, \quad
			b(\vek u, q) = \int_\Omega q \, \operatorname{div} \vek u \, \dx, \\
			& c(p,q) = \lambda^{-1} \int_\Omega p \, q \, \dx \mand 
			f_{\vek u}(\vek v) = \int_{\Gamma_\textup{right}} \vek g \cdot \vek v \, \dx.
		\end{aligned}
	\end{equation*}

	At the converged solution, the body is in equilibrium.
	At an intermediate iterate, we can interpret $r_{\vek u} = f_{\vek u} - a(\vek u,\cdot) - b(\cdot,p)$ as a net force acting on the body, measured in \si{\newton}.
	This force is attributed to a violation of the equilibrium conditions $\operatorname{div} \boldsymbol{\sigma} = - \vek{f}$ in $\Omega$, and $\boldsymbol{\sigma} \, \vek n = \Bzero$ or $\boldsymbol{\sigma} \, \vek n = \vek g$ at the non-clamping boundary parts at an intermediate iterate $(\vek u, \vek p)^{T}$.
	The second residual $r_{\vek p} = - b(\vek u,\cdot) + c(p,\cdot)$ admits an interpretation of a volume measured in \si{\meter\cubed} due to a violation of $\operatorname{div} \vek u = \lambda^{-1} p$.

	As preconditioner $\vek P = \blkdiag(\vek P_{\vek u},\vek P_{\vek p})$, we use the block diagonal matrix induced by the bilinear forms $a(\vek u,\vek v)$ and $c(p,q)$, respectively.
	Once again, we remark that this choice is compatible with the physical units of the problem.

	Similarly as for \cref{example.Stokes}, we discretize the problem using the Taylor-Hood finite element, and similar modifications due to Dirichlet displacement boundary conditions apply.
	\cref{fig:Elasticity} illustrates the convergence behavior on successively refined meshes, revealing once again different orders of magnitude for the two components of the residual.
\end{example}

\begin{figure}[htbp]
	\centering
	\includegraphics[width=0.32\textwidth]{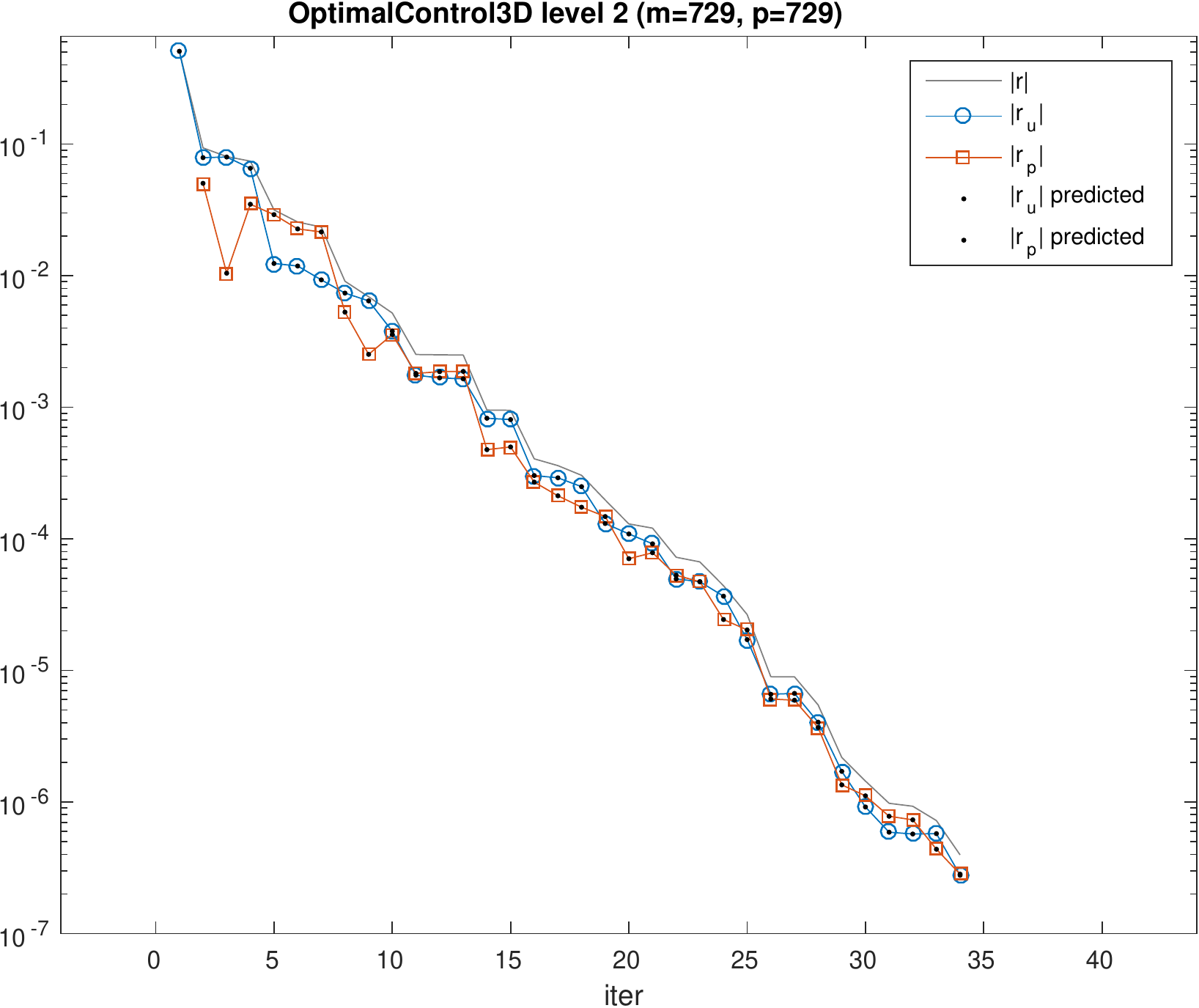}
	\includegraphics[width=0.32\textwidth]{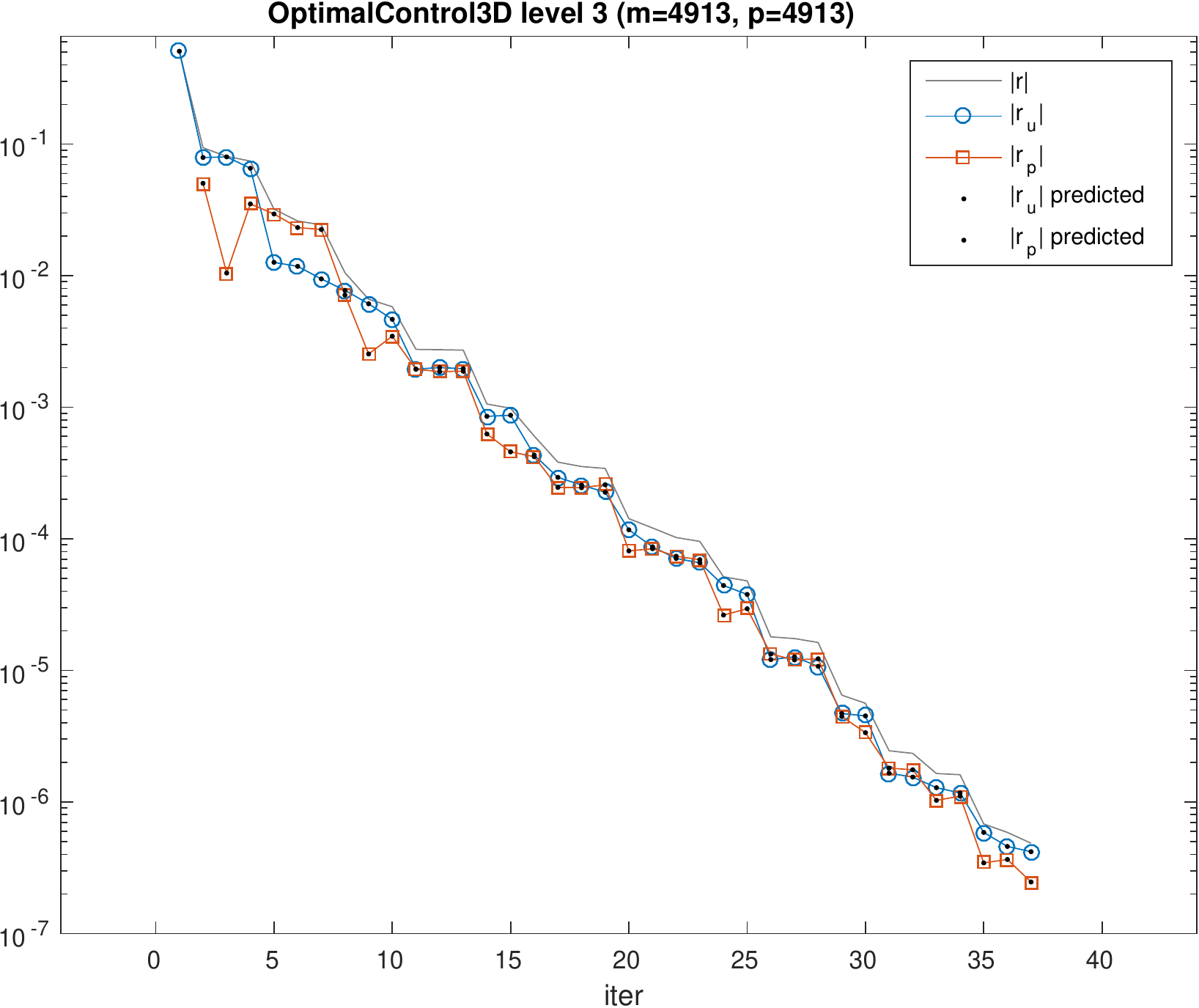}
	\includegraphics[width=0.32\textwidth]{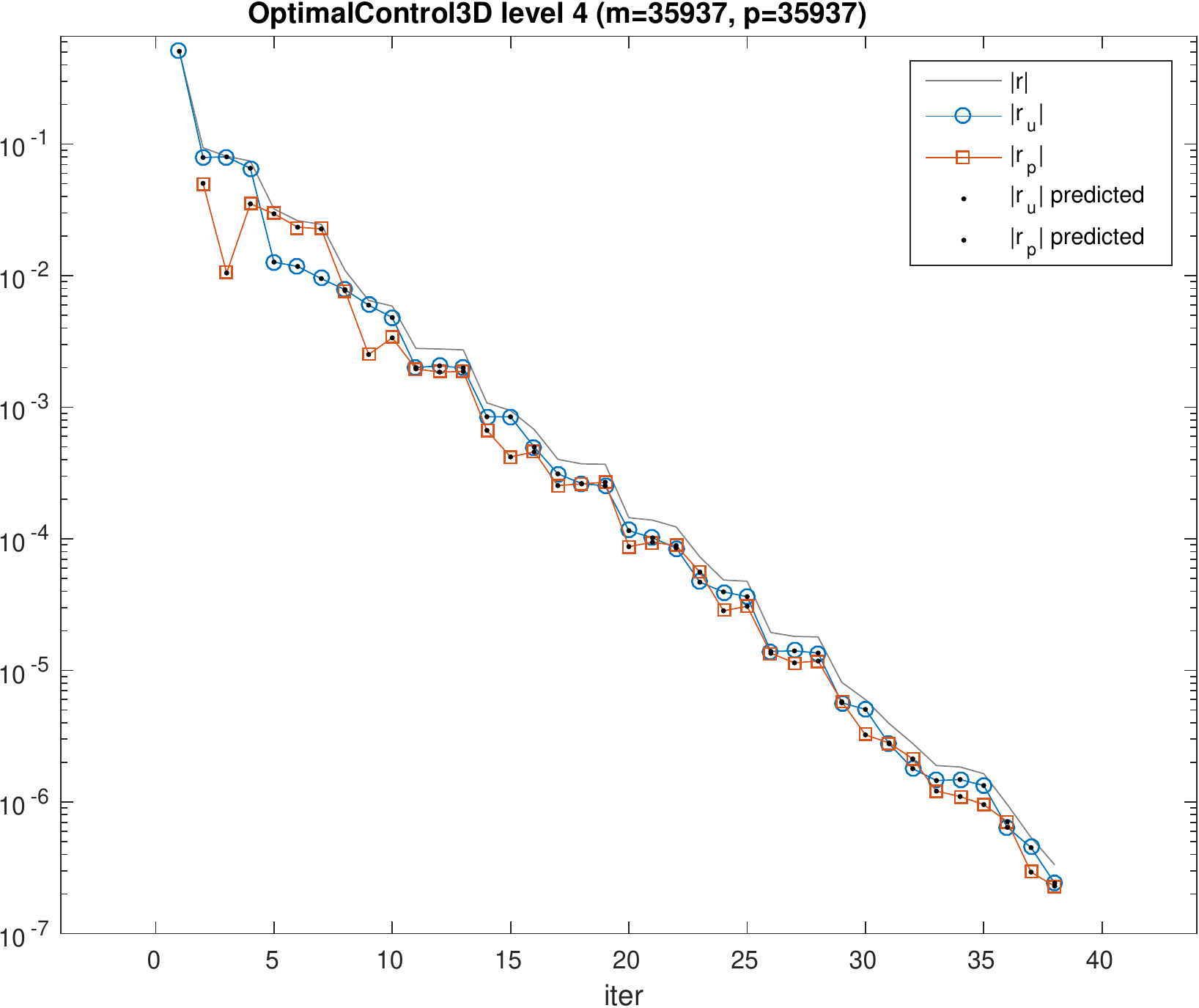}
	\caption{Convergence history of the residual subvectors for \cref{example.OptimalControl} (optimal boundary control) on a coarse grid (left) and its uniform refinements (middle and right).}
	\label{fig:OptimalControl}
\end{figure}

\begin{example}[Optimal Boundary Control]  \hfill \\
	\label{example.OptimalControl}
	Our final example is an optimal boundary control problem for the stationary heat equation on the unit cube $\Omega = (0,1) \times (0,1) \times (0,1)$ with boundary $\Gamma$.
	The problem statement is
	\begin{equation*}
		\begin{aligned}
			\text{Minimize} \quad & \frac{\alpha_1}{2} \norm{u - u_d}_{L^2(\Omega)}^2 + \frac{\alpha_2}{2} \norm{f}_{L^2(\Gamma)}^2  \\
			\text{s.t.} \quad & 
			\left\{
				\begin{aligned}
					- \kappa \, \laplace u + \delta \, u & = 0 & & \text{in } \Omega, \\
					\kappa \frac{\partial u}{\partial n} & = f & & \text{on } \Gamma.
				\end{aligned}
			\right.
		\end{aligned}
	\end{equation*}
	As data we use 
	$\alpha_1 = \SI{1}{\per\meter\cubed\per\kelvin\squared}$,
	$u_d(x,y,z) = x$~\si{\kelvin}, 
	$\alpha_2 = \SI{1e-2}{\meter\squared\per\watt\squared}$,
	heat conduction coefficient $\kappa = \SI{1}{\watt\per\meter\per\kelvin}$,
	and radiation coefficient $\delta = \SI{1}{\watt\per\meter\cubed\per\kelvin}$.

	The necessary and sufficient optimality conditions for this problem are standard; see, for instance, \cite[Chapter~2.8]{Troeltzsch2010:1}.
	When we eliminate the control function $f$ from the problem, a saddle-point system for the state $u$ and adjoint state $p$ remains, which is described by the following data:
	\begin{equation*}
		\begin{aligned}
			& a(u, v) = \alpha_1 \int_\Omega u \, v \, \dx, \quad
			b(u, q) = \kappa \int_\Omega \nabla u \cdot \nabla q \, \dx + \delta \int_\Omega u \, q \, \dx, \\
			& c(p,q) = \alpha_2^{-1} \int_\Gamma p \, q \, \dx \mand 
			f_{u}(v) = \alpha_1 \int_\Omega u_d \, v \, \dx.
		\end{aligned}
	\end{equation*}
	The underlying spaces are $V = Q = H^1(\Omega)$, and we discretize both using piecewise linear, continuous finite elements.
	Note that the elements $u, v \in V$ are measured in \si{\kelvin} while the elements $p, q \in Q$ are measured in \si{\per\watt}.

	As preconditioner $\vek P = \blkdiag(\vek P_{\vek u},\vek P_{\vek p})$, we use the block diagonal matrix induced by the bilinear forms 
	\begin{equation*}
		\begin{aligned}
			p_{\vek u}(u,v) & = \alpha_1 \kappa \, \delta^{-1} \int_\Omega \nabla u \cdot \nabla v \, \dx + \alpha_1 \int_\Omega u \, v \, \dx \\
			p_{\vek p}(p,q) & = \alpha_2^{-1} (\kappa \, \delta^{-1})^{1/2} \int_\Omega \nabla p \cdot \nabla q \, \dx + \alpha_2^{-1} (\kappa^{-1} \, \delta)^{1/2} \int_\Omega p \, q \, \dx, 
		\end{aligned}
	\end{equation*}
	respectively.
	As in our previous examples, this choice is compatible with the physical units of the problem.
	\cref{fig:OptimalControl} illustrates the convergence behavior on successively refined meshes.
	Besides the mesh independence, we observe that both residual subvector norms converge in unison in this example.
\end{example}

\section{Discussion}\label{section.discussion}
In this paper we developed a modified implementation of MINRES.
When applied to saddle-point systems, the new implementation allows us to monitor the norms of the subvectors $\bignorm{\vek r_{\vek u}^{(j)}}_{\vek P_{\vek u}^{-1}}$ and $\bignorm{\vek r^{(j)}_{\vek p}}_{\vek P_{\vek p}^{-1}}$ individually, while conventional implementations keep track of only the total residual norm $\bignorm{\vek r^{(j)}}_{\vek P^{-1}}$.
It should be obvious how \cref{alg.prec.minres-modified} generalizes to systems with more than two residual subvectors and block-diagonal preconditioners structured accordingly.
The price to pay to monitor the subvector norms is the storage of one additional vector for components $\vek m_{j,\vek u}$ and $\vek m_{j,\vek p}$ compared to the implementation of MINRES given in \cite[Algorithm~4.1]{ESW-Book.2014} 
as well as some additional inner product calculations.  It should be noted that if the cost of applying the operator
and preconditioner are cheap enough, it may be preferable to simply construct the full residual and calculate preconditioned
subvector norms explicitly.  However, in the case that application of these operators has significant costs (particularly in the
case that they are available only as procedures and never actually constructed) the method presented in this paper has
important advantages.

While we developed the details in finite dimensions using matrices and vectors, our approach directly transfers to linear saddle-point systems in a Hilbert space setting using linear operators and linear forms, as described in \cite{G.Herzog.S.2014}.

Being able to differentiate between the contributions to the total residual offers new opportunities for the design of iterative algorithms for nonlinear problems, which require inexact solves of \eqref{eqn.ABupf} as their main ingredient.
We envision for instance solvers for equality constrained nonlinear optimization problems which may now assign individual stopping criteria for the residuals representing optimality and feasibility, respectively.
The design of such an algorithm is, however, beyond the scope of this work.
Similarly, for the Stokes \cref{example.Stokes}, the user may now assign individual stopping criteria for the fulfillment of the balance of forces (first residual) and the conservation of mass (second residual).

A topic of future research could be to study the decay properties of the subvector norms, rather than study the decay of the total residual norm; see for instance \cite[Section~4.2.4]{ESW-Book.2014}.
We expect that such an analysis will be more involved but it may shed light on how the relative scaling of the preconditioners blocks affects the convergence of the residual subvectors.

\section*{Acknowledgment}
The authors would like to thank David Silvester and Valeria Simoncini for their comments on an earlier version of the manuscript.

\bibliographystyle{siamplain}
\bibliography{master}
\end{document}